\documentclass[leqno,draft,11pt,a4paper]{article}
\usepackage{amssymb,amsmath,theorem}
\usepackage[includeheadfoot,margin=1in]{geometry}

\newcommand\mclap{\mathpalette\doclap}
\newcommand\doclap[2]{\hbox to0pt{\hss$#1#2$\hss}}

\newcommand\et{\mathrel\&}
\newcommand\eq{\leftrightarrow}
\newcommand\ET{\bigwedge}
\newcommand\model{\vDash}
\newcommand\nmodel{\nvDash}
\newcommand\lang{\mathcal L}
\newcommand\langor{\lang_{\mathrm{OR}}}
\newcommand\langp{\langor\cup\{P_2\}}
\newcommand\lange{\langor\cup\{2^x\}}

\newcommand\fii{\varphi}

\newcommand\p[1]{\langle#1\rangle}
\newcommand\abs[1]{\lvert#1\rvert}
\newcommand\dlh[1]{\lVert#1\rVert}
\newcommand\bez{\smallsetminus}
\newcommand\sset{\subseteq}
\newcommand\nsset{\nsubseteq}
\newcommand\Sset{\supseteq}
\newcommand\res{\mathbin\restriction}

\newcommand\delim[4]{\ifx X#3X\left#1#4\right#2\else\csname#3l\endcsname#1#4\csname#3r\endcsname#2\fi}
\newcommand\fl[1]{\lfloor#1\rfloor}
\newcommand\Fl[2][]{\delim\lfloor\rfloor{#1}{#2}}
\newcommand\cl[1]{\lceil#1\rceil}
\newcommand\CL[2][]{\delim\lceil\rceil{#1}{#2}}
\newcommand\tdive[2]{\fl{#1/2^{#2}}}

\DeclareMathOperator\im{im}
\DeclareMathOperator\Th{Th}
\DeclareMathOperator\lcm{lcm}

\newcommand\thry[1]{\mathsf{#1}}
\newcommand\delz{\Delta_0}
\newcommand\io{\thry{IOpen}}
\newcommand\idz{\thry{I\delz}}
\newcommand\idze{\idz+\thry{EXP}}
\newcommand\dicr{\thry{\Delta^b_1\text-CR}}
\newcommand\vtc{\thry{VTC^0}}
\newcommand\PA{\thry{PA}}
\newcommand\teip{\thry{TEIP}}
\newcommand\teipp{\teip_{\!P_2}}
\newcommand\teipe{\teip_{2^x}}
\newcommand\teippow{\teip_{\powii}}
\newcommand\powii{\mathrm{Pow}_2}
\newcommand\pwin{\mathrm{PWin}}
\newcommand\powg{\mathrm{PowG}}

\newcommand\N{\mathbb N}
\newcommand\Q{\mathbb Q}
\newcommand\Z{\mathbb Z}
\newcommand\R{\mathbb R}
\newcommand\sM{\mathfrak M}
\newcommand\sR{\mathfrak R}
\newcommand\sF{\mathfrak F}

\newcommand\ob[1]{\overline{#1}}
\newcommand\txto{${}\to{}$}

\newcommand\bme{\hskip.75em\relax}
\newcommand\noproof{\leavevmode\unskip\bme\vadjust{}\nobreak\hfill$\Box$\par}
\newenvironment{Pf}
  {\par\noindent\textit{Proof:}\bme\ignorespaces}
  {\noproof\pagebreak[2]\vskip\medskipamount\ignorespacesafterend}
\newcommand\qedhere{\noproof\let\noproof\relax}

\theoremstyle{plain}
\newtheorem{Thm}{Theorem}[section]
\newtheorem{Prop}[Thm]{Proposition}
\newtheorem{Cor}[Thm]{Corollary}
\newtheorem{Lem}[Thm]{Lemma}
\newtheorem{Obs}[Thm]{Observation}
\newtheorem{Que}[Thm]{Question}
\theorembodyfont\upshape
\newtheorem{Def}[Thm]{Definition}
\newtheorem{Rem}[Thm]{Remark}
\newtheorem{Exm}[Thm]{Example}

\author{Emil Je\v r\'abek%
   \thanks{Supported by the Czech Academy of Sciences (RVO 67985840) and GA \v CR project 23-04825S.}\\[\medskipamount]
Institute of Mathematics, Czech Academy of Sciences\\
\small \v Zitn\'a 25,
110\:00 Praha 1,
Czech Republic,
email: \texttt{jerabek@math.cas.cz}
}

\title{On the theory of exponential integer parts}

\begin{document}
\maketitle

\begin{abstract}
We axiomatize the first-order theories of exponential integer parts of real-closed exponential fields in a language
with $2^x$, in a language with a predicate for powers of~$2$, and in the basic language of ordered rings. In
particular, the last theory extends $\io$ by sentences expressing the existence of winning strategies in a certain game
on integers; we show that it is a proper extension of $\io$, and give upper and lower bounds on the required number of
rounds needed to win the game.
\end{abstract}

\section{Introduction}\label{sec:introduction}

A classical result of Shepherdson~\cite{sheph} characterizes models of the arithmetical theory $\io$ (induction for
quantifier-free formulas in the language $\langor=\p{0,1,+,\cdot,<}$) as exactly those that are (nonnegative parts of)
integer parts of real-closed fields. Here, an integer part (IP) of an ordered ring~$R$ is a discrete subring $I\sset R$
such that every element of~$R$ is within distance~$1$ from an element of~$I$. An analogue for exponential ordered
fields $\p{R,\exp}$ (with $\exp(1)=2$) was introduced by Ressayre~\cite{ress:eip}: an exponential integer part (EIP)
of~$R$ is an IP $I\sset R$ such that $I_{\ge0}$ is closed under $\exp$. (We will find it more convenient to call the
nonnegative part $I_{\ge0}$ the EIP of~$R$ rather than $I$ itself, and we usually denote ${\exp}\res I_{\ge0}$ as
$2^x$.) We are interested in the question what models of $\io$ are EIP of real-closed exponential fields (RCEF), and in
particular, what is the first-order theory of such structures.

The question whether the theory of EIP of RCEF in~$\langor$ coincides with $\io$ was raised by Je\v
r\'abek~\cite{ej:vtceip}; he provided an upper bound on the theory, proving that every countable model of the weak
arithmetic $\vtc$ (or rather, the equivalent one-sorted arithmetical theory $\dicr$) is an EIP of a RCEF, despite the
fact that the ``natural'' integer exponentiation function in this theory is only defined for small integers. (This
considerably improved the base-line bound given by the more-or-less obvious fact%
\footnote{The formalization of integer powering and other elementary recursive functions in $\idze$ implies that it can
define well-behaved rational approximations of $\exp(x)$ for rational~$x$, which implies that the completion of a
fraction field of a model of $\idze$---which is a RCF---carries an exponential. Alternatively, one might restrict the
completion to limits of elementary recursive Cauchy sequences, though this requires formalizing in $\idze$ basic facts
about polynomial root finding to ensure the resulting exponential field is real-closed. Surprisingly, this seems not to
have been explicitly pointed out in the literature until recently; it is stated in~\cite{ej:vtceip} as a consequence of
the results from~\cite{ej:vtcanal} on formalization of $\exp$ in $\vtc$ (which is an overkill for $\idze$). Earlier,
Carl and Krapp~\cite{carl-krapp} published a proof of the much weaker claim that models of $\PA$ are EIP of RCEF;
strengthening this in another direction, they proved that every model of $\Th(\N)$ is an EIP of a model of
$\Th(\R,\exp)$.}
that every model of $\idze$ is an EIP of a RCEF.)

Extensions of Shepherdson's theorem to RCEF were studied previously by Boughattas and Ressayre~\cite{bou-ress:eip} and
Kovalyov~\cite{kov:eip}. Their work differs from ours in two main respects. First, they approach Shepherdson's
characterization from the other side, focusing on problems such as: what additional axioms must be added to RCEF to
ensure that their EIP are models of such and such theory (e.g., open induction in a language with exponentiation).
Second, they mostly study IP closed under the binary powering operation $x^y=\exp(y\log x)$ in a language including
$x^y$: in this case, the ambient exponential field can be canonically reconstructed from the structure of the IP (using
the integer $x^y$ operation, we can define rational approximations of $\exp(x)$ for rational~$x$, which yields an
exponential function on the completion of the fraction field). Such a direct construction seems impossible if we
have only $2^x$ instead of $x^y$ in the language, let alone when we work only with the basic language~$\langor$; thus,
we will instead rely on model-theoretic tools such as the joint consistency theorem and recursively saturated models.

Our main contribution is an axiomatization of the theories of EIP of RCEF in $\lange$, in $\langp$ (where $P_2$ is a
predicate for the set of powers of~$2$), and in~$\langor$, denoted $\teipe$, $\teipp$, and $\teip$ (respectively). The
first two theories are extensions of $\io$ by finitely many axioms expressing basic algebraic properties of $2^x$
and~$P_2$. The most important theory, $\teip$, is more involved: it extends $\io$ with an infinite sequence of
sentences expressing that a certain game on positive integers (designed so that playing powers of~$2$ is a winning
strategy) is a win for the second player. We note that there is a general result on axiomatizing conservative fragments
of given theories by means of game sentences of similar kind due to Svenonius~\cite{sve:extra}, which is instrumental
in the argument that countable recursively saturated models are resplendent (Barwise and Schlipf~\cite{bar-sch}).
However, in contrast to the rather opaque game considered by Svenonius, mimicking the Henkin completion procedure, our
game on integers has simple and transparent rules, which makes it amenable to combinatorial analysis.

We show that $\teip$ is a proper extension of $\io$. We leave open the problem whether $\teip$ is finitely
axiomatizable over $\io$, but as a partial progress, we prove that formulas obtained by stripping the outermost pair of
quantifiers from each axiom of $\teip$ form a strict hierarchy (even over the true arithmetic~$\Th(\N)$); this amounts
to the fact that if we play our integer game starting with arbitrary numbers that are not powers of~$2$, the first
player needs an unbounded number of rounds to win. To this end, we analyze the game, proving upper and lower bounds on
the number of rounds needed to win that are tight for a sizeable set of initial integers.

There is a natural interpretation of $\langp$ in arithmetic where we put $P_2(x)$ iff $x$ has no nontrivial odd
divisor (``$x$ is oddless''). We briefly discuss what theories of arithmetic prove $\teipp$ under this interpretation
(and hence include $\teip$): in particular, this holds for $\thry{IE_2}$. On the other hand, not even $\Th(\N)+\teipp$
can prove that $P_2(x)$ implies $x$ is oddless: it is consistent that an element of $P_2$ is divisible by~$3$. In other
words, even for strong theories of arithmetic, expansions to models of $\teipp$ are not unique.

The paper is organized as follows. We review the preliminaries in Section~\ref{sec:preliminaries}. We compute the
theories of EIP of RCEF in $\lange$ and~$\langp$ in Section~\ref{sec:teipe-teipp}. In Section~\ref{sec:teip}, we
introduce the $\powg$ game and determine the theory $\teip$ of EIP or RCEF in~$\langor$.
Section~\ref{sec:analysis-powg} is devoted to an analysis of winning strategies in~$\powg$. We discuss the oddless
interpretation of~$P_2$ in Section~\ref{sec:oddless}, and we end with some concluding remarks and open problems in
Section~\ref{sec:conclusion}.

\section{Preliminaries}\label{sec:preliminaries}
Let $\langor=\p{0,1,+,\cdot,{<}}$. An \emph{ordered ring} is an $\langor$-structure $\sR=\p{R,0,1,+,\cdot,{<}}$ such
that $\p{R,0,1,+,\cdot}$ is a commutative ring, and $<$ is a (strict) total order on~$R$ compatible with $+$
and~$\cdot$ (i.e., $x\le y$ implies $x+z\le y+z$, and, if $z\ge0$, also $xz\le yz$). We call $\sR$ \emph{discrete} if
there is no element strictly between $0$ and~$1$. An \emph{ordered field} is an ordered ring that is a field. An
ordered field~$\sR=\p{R,\dots}$ is a \emph{real-closed field} (\emph{RCF}) if it has no proper algebraic extension to an
ordered field, or equivalently, if every $a\in R_{>0}$ has a square root in~$\sR$, and every polynomial $f\in\sR[x]$ of
odd degree has a root in~$\sR$, where $R_{>0}$ denotes $\{a\in R:a>0\}$.

An \emph{integer part} (\emph{IP}) of an ordered ring~$\sR=\p{R,\dots}$ is a discrete subring $I\sset R$ (considered as an
$\langor$-substructure) such that for every $a\in R$, there is $z\in I$ such that $z\le a<z+1$.

The \emph{nonnegative part} of an ordered ring $\sR=\p{R,\dots}$ is its substructure $\sR_{\ge0}$ with domain
$R_{\ge0}=\{a\in R:a\ge0\}$. The theory of nonnegative parts of discrete ordered rings is denoted $\PA^-$; it is an
extension of Robinson's arithmetic~$\thry Q$. Every $\sM=\p{M,\dots}\model\PA^-$ has a unique (up to isomorphism)
extension to a discrete ordered ring $\sM_\pm=\p{M_\pm,\dots}$ such that $(\sM_\pm)_{\ge0}=\sM$ and
$M_\pm=\{a,-a:a\in M\}$, which is called the \emph{extension of $\sM$ with negatives}. The extension of $\PA^-$ (or
equivalently, $\thry Q$) with the \emph{induction axioms}
\begin{equation}\label{eq:ind}\tag{$\fii$-IND}
\forall\vec y\:\bigl(\fii(0,\vec y)\land\forall x\:\bigl(\fii(x,\vec y)\to\fii(x+1,\vec y)\bigr)
   \to\forall x\,\fii(x,\vec y)
\end{equation}
for all open (= quantifier-free) $\langor$-formulas~$\fii$ is denoted $\io$.
\begin{Thm}[Shepherdson \cite{sheph}]\label{thm:shep}
An $\langor$-structure $\sM$ is an IP of a RCF if and only if $\sM_{\ge0}\model\io$.
\noproof\end{Thm}
Note that a priori there is no reason for the class of integer parts of RCF to be elementary; indeed, this fails for
our case of interest (EIP of RCEF), as we will see.

$\sR=\p{R,0,1,+,\cdot,{<},{\exp}}$ is an \emph{(ordered) exponential field} if $\p{R,0,1,+,\cdot,{<}}$ is an ordered
field, and $\exp$ is an ordered group isomorphism $\exp\colon\p{R,+,0,{<}}\to\p{R_{>0},\cdot,1,{<}}$. Following
Ressayre~\cite{ress:eip}, a \emph{real-closed exponential field} (\emph{RCEF}) is an exponential field $\sR=\p{R,\dots}$ which
is real-closed and satisfies $\exp(1)=2$; if $\exp(x)>x$ for all $x\in R$, we say that it satisfies the \emph{growth
axiom}%
\footnote{Ressayre includes this in the definition of an exponential field, and actually formulates it as
``$\exp(x)>x^n$ for all $x$ somewhat larger than $n$'', where $n$ presumably refers to standard natural numbers. This
follows from our GA, since $\exp(x)=\exp(x/2n)^{2n}>(x/2n)^{2n}\ge x^n$ as long as, say, $x\ge(2n)^2$. On the other
hand, it is easy to see that if there is $m\in\N$ such that $\exp(x)>x$ holds for all $x\ge m$, then it holds for all
$x\in R$, thus our axiom is equivalent to Ressayre's formulation.}
(\emph{GA}). If $I$ is an IP of $\sR$ such that $I_{\ge0}$ is closed under~$\exp$, we call $I_{\ge0}$ an \emph{exponential
integer part} (\emph{EIP}) of~$\sR$. (We define $I_{\ge0}$, rather than $I$ itself, to be an EIP, since we intend to
axiomatize first-order theories of EIP as extensions of $\io$, and compare them with other theories of arithmetic such
as $\thry{IE}_k$, which are formulated such that all elements are nonnegative.) We consider an EIP $I_{\ge0}$ to be not
just a set, but an $\langor$-substructure of~$\sR$, and we also consider it in some expanded languages: $\lange$, by
inheriting the function $2^x={\exp}\res I_{\ge0}$ from~$\sR$, and $\langp$, where the unary predicate $P_2$ is
interpreted as the image of $2^x\colon I_{\ge0}\to I_{>0}$.

\emph{Presburger arithmetic} is the complete theory of the structure $\p{\N,0,1,+,<}$. Models of Presburger
arithmetic are exactly the nonnegative parts of \emph{$\Z$-groups}, which are discrete ordered abelian groups
$\p{Z,0,+,<}$ with a least positive element~$1$ such that $Z/\Z$ is divisible, where we identify $\Z$ with the subgroup
of~$Z$ generated by~$1$. There is an (easily proved) baby version of Theorem~\ref{thm:shep}: $\Z$-groups are exactly the IP
of divisible ordered abelian groups (where an IP of an ordered group is defined analogously to rings, but without
multiplication).

In theories extending $\PA^-$, \emph{existential bounded quantifiers} $\exists x\le t\,\fii(x,\dots)$ (where $t$ is a
term that does not contain~$x$) are defined as shorthands for $\exists x\,(x\le t\land\fii(x,\dots))$, and
\emph{universal bounded quantifiers} $\forall x\le t\,\fii(x,\dots)$ are shorthands for
$\forall x\,(x\le t\to\fii(x,\dots))$. A \emph{bounded formula} is one that only uses bounded quantifiers. The set of
all bounded $\langor$-formulas is denoted $\delz$. An $\langor$-formula is $E_k$ (resp., $U_k$) if it can be written
with $k$~alternating (possibly empty) blocks of bounded quantifiers followed by a quantifier-free formula, with the
first block being existential (resp., universal). If $\Gamma$ is a formula class such as $\delz$ or~$E_k$,
$\thry I\Gamma$ denotes the theory axiomatized by $\PA^-$ (or just $\thry Q$) and \eqref{eq:ind} for formulas
$\fii\in\Gamma$ (thus, $\thry{IE}_0=\io$).

We define the divisibility predicate $x\mid y$ as $\exists z\,xz=y$ (thus all elements divide~$0$). Over $\PA^-$, the
existential quantifier can be bounded by $z\le y$, thus $x\mid y$ is an $E_1$~formula; it is equivalent to the
$U_1$~formula $\forall q\le y\,\forall r<x\,(y=qx+r\to r=0)$ over $\io$.

The theory $\dicr$ of Johannsen and Pollett~\cite{joh-pol:d1cr} is a weak theory of bounded arithmetic in the style of
Buss's theories (cf.\ \cite[\S V.4]{hp}) that corresponds to the complexity class $\mathrm{TC}^0$. It is
bi-interpretable (RSUV-isomorphic) to the more commonly used two-sorted Zambella-style theory $\vtc$ (see
\cite{cook-ngu}), but since our interest lies in embedding the universe of the theory with its $\langor$-structure as
EIP in other structures, it is more natural to consider the one-sorted version of the theory. It was proved in Je\v
r\'abek~\cite{ej:vtceip} that every countable model of $\dicr$ is an EIP of a RCEF satisfying GA (despite the fact that
the natural exponentiation function in $\dicr$ is only defined on an initial segment of small integers). Proper
definitions of $\dicr$ and $\vtc$ as well as more context can be found in the references above; readers unfamiliar with
these theories may safely skip the few places where they are mentioned below.

We use $\log x$ to denote the base-$2$ logarithm of~$x$, with the convention that $\log x=0$ for $x\le1$ (i.e., it is
really $\max\{0,\log_2x\}$). We denote the natural logarithm by $\ln x$, and general base-$b$ logarithm by $\log_bx$.

We will also need two tools from model theory. The first is Robinson's joint consistency theorem (see e.g.\ Hodges
\cite[Cor.~9.5.8]{hodges:mod-th}):
\begin{Thm}\label{thm:jct}
Let $T$ be a complete $\lang$-theory, and for $i=0,1$, let $T_i\Sset T$ be a consistent $\lang_i$-theory, where
$\lang_0\cap\lang_1=\lang$. Then $T_0\cup T_1$ is consistent.
\noproof\end{Thm}

Recursive saturation was introduced by Barwise and Schlipf~\cite{bar-sch}. Let $\sM=\p{M,\dots}$ be a
structure in a finite language~$\lang$. If $\vec a\in M$ and $\Gamma(x,\vec y)$ is a recursive set of $\lang$-formulas,
then $\Gamma(x,\vec a)$ is a \emph{recursive type} of~$\sM$, which is \emph{finitely satisfiable} if
$\sM\model\exists x\,\ET_{\fii\in\Gamma'}\fii(x,\vec a)$ for each finite $\Gamma'\sset\Gamma$, and
\emph{realized} by $c\in M$ if $\sM\model\Gamma(c,\vec a)$. Then $\sM$ is \emph{recursively saturated} if
every finitely satisfiable recursive type of~$\sM$ is realized in~$\sM$.

For two structures $\sM$ and $\sM'$, $\p{\sM,\sM'}$ denotes the two-sorted structure that combines them, which can be
represented by a one-sorted structure in the usual way (make the languages relational, let the domain be the
disjoint union of $M$ and $M'$, and include a unary predicate that separates $M$ from~$M'$). $\sM$ and~$\sM'$ are
\emph{jointly recursively saturated} if $\p{\sM,\sM'}$ is recursively saturated. E.g., two reducts (possibly
in a renamed language) of the same recursively saturated structure, or more generally, two structures interpretable in
one recursively saturated structure, are jointly recursively saturated.

We will use the basic existence and uniqueness results on recursively saturated structures:
\begin{Lem}[{{Barwise and Schlipf~\cite[\S1.4]{bar-sch}}}]\label{lem:rec-sat}
\
\begin{enumerate}
\item Every structure has a recursively saturated elementary extension of the same cardinality.
\item If $\sM$ and~$\sM'$ are elementarily equivalent jointly recursively saturated countable structures, then
$\sM\simeq\sM'$.\qedhere
\end{enumerate}
\end{Lem}

We note that Theorem~\ref{thm:jct} follows almost immediately from Lemma~\ref{lem:rec-sat} as well.

\section{Exponential integer parts in a language with $2^x$ or $P_2$}\label{sec:teipe-teipp}

We start by axiomatizing the theory of EIP of RCEF in a language with~$2^x$, which is fairly straightforward.

\begin{Def}\label{def:teipe}
$\teipe$ is a theory in the language $\lange$ extending $\io$ by the axioms
\begin{align}
\label{eq:1}\tag{$2^x$-IP}&x>0\to\exists y\:x<2^y\le2x,\\
\label{eq:2}\tag{$2^x$-Mul}&2^{x+y}=2^x2^y,\\
\label{eq:3}\tag{$2^x$-Pos}&2^x>0.
\intertext{$\teipe^+$ is defined similarly, but with axiom}
\tag{$2^x$-GA}\label{eq:3+}&2^x>x
\end{align}
in place of~\eqref{eq:3}.
\end{Def}

By doubling/halving $x$ (which corresponds to shifting $y$ by~$1$), \eqref{eq:1} is equivalent to
\[x>0\to\exists y\:2^y\le x<2\cdot2^y,\]
which would match more closely the axioms of $\teipp$ and $\teip$ that will be given further on, but the version here
looks more visually pleasing.

\begin{Thm}\label{thm:teipe}
The first-order theory of EIP of RCEF in $\lange$ is $\teipe$. 
The first-order theory of EIP of RCEF satisfying GA in $\lange$ is $\teipe^+$. 
\end{Thm}
\begin{Pf}
It is clear that any EIP of a RCEF satisfies the given axioms. Conversely, assume that
$\sM=\p{M,0,1,+,\cdot,<,2^x}\model\teipe$. Since $\sM\model\io$, $\sM_\pm$ is an IP of a RCF~$\sR$ by
Theorem~\ref{thm:shep}. There exists an elementary extension $\sR^*=\p{R^*,M^*,0,1,+,\cdot,<,2^x}$ of
$\p{R,M,0,1,+,\cdot,<,2^x}$ that expands to a RCEF $\p{\sR^*,\exp}$ by Theorem~\ref{thm:jct} (applied with $T_0$ being the
elementary diagram of $\p{R,M,0,1,+,\cdot,<,2^x}$ and $T_1$ the theory of RCEF, with common language~$\langor$), using
the completeness of the theory RCF\@. Let $M^*_\pm=M^*\cup\{-a:a\in M\}$, and extend $2^x\colon M^*\to M^*_{>0}$ to a
function $2^x\colon M^*_\pm\to R^*_{>0}$ by $2^{-x}=(2^x)^{-1}$. Applying \eqref{eq:1} with $x=1$, there exists
$y\in M^*$ such that $2^y=2$; depending on the parity of~$y$, \eqref{eq:2} implies $2=(2^{y/2})^2$ (which is
impossible) or $2=2^1\cdot(2^{\fl{y/2}})^2$, thus $2^1=2$. Then using \eqref{eq:2} and~\eqref{eq:3}, $2^x$ is strictly
increasing, hence it is an ordered group embedding $\p{M^*_\pm,0,+,<}\to\p{R^*_{>0},1,\cdot,<}$.

Putting $B=\{x\in R^*:\exists n\in\N\,\abs x\le n\}$, we define a new exponential $\ob\exp\colon R^*\to R^*_{>0}$ by
\[\ob\exp(a+r)=2^a\exp(r),\qquad a\in M^*_\pm,r\in B.\]
To see that this is well defined, if $a+r=a'+r'$ with $a,a'\in M^*_\pm$ and $r,r'\in B$, then $n=a-a'=r'-r\in B\cap
M^*_\pm=\Z$, hence both $2^{a-a'}$ and $\exp(r'-r)$ coincide with the usual value of $2^n$, which implies
$2^a\exp(r)=2^{a'}\exp(r')$. The function $\ob\exp$ is defined on all of~$R^*$ as $M^*_\pm$ is an IP of~$R^*$. 

It follows easily that $\ob\exp$ is a homomorphism $\p{R^*,0,+}\to\p{R^*_{>0},1,\cdot}$ using the corresponding
properties of $2^x$ and $\exp$, and $\ob\exp(x)>1$ for $x>0$, thus $\ob\exp$ is strictly increasing. It is also
surjective: if $y\in R^*_{>0}$, \eqref{eq:1} implies that there is $a\in M^*_\pm$ such that $2^{-a}y\in[1,2]$,
thus using the surjectivity and monotonicity of $\exp$, $2^{-a}y=\exp(x)$ for some $x\in[0,1]$, whence
$y=\ob\exp(a+x)$. That is, $\p{R^*,0,1,+,\cdot,<,\ob\exp}$ is an RCEF, and $\p{M^*_\pm,0,1,+,\cdot,<,2^x}$ is its EIP.

If $\sM$ and $\sM^*$ additionally satisfy~\eqref{eq:3+}, then $\ob\exp(x)>x$ for all $x$: this holds trivially if
$x\le0$; otherwise, we can write $x=a+r$ with $a\in M$ and $r\in[0,1)$, thus $2^x\ge2^a\ge a+1>x$.
\end{Pf}

Next, we move to a language that only has a predicate $P_2$ for the image of~$2^x$ rather than $2^x$ itself. It turns
out that the resulting theory is the same irrespective of whether we demand the RCEF to satisfy the growth axiom; but
whereas in absence of GA, the proof is still straightforward, the GA case is considerably more complicated.

\begin{Def}\label{def:teipp}
$\teipp$ is a theory in the language $\langp$ extending $\io$ by the axioms
\begin{align}
\label{eq:4}\tag{$P_2$-IP}&x>0\to\exists u\:(P_2(u)\land u\le x<2u),\\
\label{eq:5}\tag{$P_2$-Div}&P_2(u)\land P_2(v)\land u\le v\to\exists w\:(P_2(w)\land uw=v).
\intertext{$\teipp'$ is a theory in the language $\langp$ extending $\io$ by the axioms}
\label{eq:6}\tag{$P_2$-IP!}&x>0\to\exists!u\:(P_2(u)\land u\le x<2u),\\
\label{eq:7}\tag{$P_2$-Pos}&\neg P_2(0),\\
\label{eq:8}\tag{$P_2$-Mul}&P_2(u)\land P_2(v)\to P_2(uv).
\end{align}
\end{Def}

\begin{Lem}\label{lem:teipp}
$\teipp$ is equivalent to $\teipp'$, and it proves $P_2(1)$ and $P_2(2)$.
\end{Lem}
\begin{Pf}

$\teipp\vdash\teipp'$: The existence part of~\eqref{eq:6} is just~\eqref{eq:4}; for uniqueness, if $u$ and~$u'$ satisfy
the conclusion, and, say, $u\le u'$, then $u\mid u'$ by~\eqref{eq:5}, while $u'<2u$. Thus, the only possibility
is $u=u'$.

Applying \eqref{eq:4} with $x=1$, we see that $P_2(1)$. Then \eqref{eq:5} gives $\neg P_2(0)$ as $0\nmid1$.

Assume $P_2(u)$ and $P_2(v)$. If $u=1$ or $v=1$, then $P_2(uv)$ holds trivially, hence we may also assume $u,v\ge2$.
By~\eqref{eq:4}, there is $w$ such that $P_2(w)$ and $w\le uv<2w$. Since $2u\le uv$, we have $u<w$, hence \eqref{eq:5}
implies $u\mid w$ and $P_2(w/u)$. Moreover, $w/u\le v<2w/u$, thus $w/u=v$ by the uniqueness part of~\eqref{eq:6}, i.e.,
$P_2(uv)$.

$\teipp'\vdash\teipp$: \eqref{eq:4} follows from~\eqref{eq:6}. Assume that $P_2(u)$, $P_2(v)$, and $u\le v$. We have
$u>0$ by \eqref{eq:7}, whence $\io$ implies the existence of $x>0$ such that $ux\le v<u(x+1)$. 
By~\eqref{eq:6}, there is $w$ such that $P_2(w)$ and $w\le x<2w$, i.e., $uw\le v<2uw$. Then $P_2(uw)$ by~\eqref{eq:8},
hence $uw=v$ by the uniqueness part of~\eqref{eq:6}.

We have already seen that $\teipp\vdash P_2(1)$. Likewise, an application of~\eqref{eq:4} with $x=2$ gives $P_2(2)$.
\end{Pf}

\begin{Lem}\label{lem:teipp-pres}
If $\sM=\p{M,0,1,+,\cdot,<,P_2}\model\teipp$, then $\p{P_2,1,2,\cdot,<}$ is a model of Presburger arithmetic.
\end{Lem}
\begin{Pf}
Let $\sR=\p{R,\dots}$ be a RCF such that $\sM_\pm$ is its IP using Theorem~\ref{thm:shep}. Let
$P_2^{\pm1}=\{u,u^{-1}:u\in P_2\}\sset R_{>0}$. Using Lemma~\ref{lem:teipp}, $\p{P_2^{\pm1},1,\cdot,<}$ is a discrete
ordered abelian group with a least positive element~$2$, and it is an IP of the divisible ordered group
$\p{R_{>0},1,\cdot,<}$. Thus, it is a $\Z$-group, and its ``nonnegative'' part $\p{P_2,1,2,\cdot,<}$ is a model of
Presburger.
\end{Pf}

For the construction of $\teip$ in the next section, it will be convenient to consider yet another axiomatization of
$\teipp$ that may look less intuitive, but has the advantage that it only involves one positive occurrence of~$P_2$:
\begin{Def}\label{def:teippp}
$\teipp''$ is a theory in the language $\langp$ extending $\io$ by the axioms \eqref{eq:4}, \eqref{eq:7}, and
\begin{equation}
\label{eq:11}\tag{$P_2$-Univ}P_2(u)\land P_2(v)\land P_2(w)\to\neg(uv<w<2uv).
\end{equation}
\end{Def}

\begin{Lem}\label{lem:teippp}
$\teipp'$ is equivalent to $\teipp''$.
\end{Lem}
\begin{Pf}

$\teipp'\vdash\teipp''$: For~\eqref{eq:11}, if $P_2(u)$, $P_2(v)$, and~$P_2(w)$, then $P_2(uv)$ by~\eqref{eq:8}, hence
the uniqueness part of~\eqref{eq:6} precludes $uv<w<2uv$.

$\teipp''\vdash\teipp'$: First, \eqref{eq:4} implies $P_2(1)$, hence \eqref{eq:11} gives $P_2(u)\land
P_2(u')\to\neg(u<u'<2u)$, which is the uniqueness part of~\eqref{eq:6}; thus, we have \eqref{eq:6} and~\eqref{eq:7}.
For~\eqref{eq:8}, assuming $P_2(u)$ and $P_2(v)$, there is $w$ such that $P_2(w)$ and $uv\le w<2uv$ by~\eqref{eq:4}; we
must have $w=uv$ by~\eqref{eq:11}.
\end{Pf}

\begin{Thm}\label{thm:teipp}
The first-order theory of EIP of RCEF in $\langp$ is $\teipp$. 
\end{Thm}
\begin{Pf}
In view of Theorem~\ref{thm:teipe}, it suffices to show that $\teipe$ is a conservative extension of $\teipp$ when $P_2(u)$
is interpreted as $\exists x\,u=2^x$. Clearly, $\teipe$ proves $\teipp$.

On the other hand, let $\sM=\p{M,0,1,+,\cdot,<,P_2}\model\teipp$, which we may assume to be countable and recursively
saturated. Using $\io$ and Lemma~\ref{lem:teipp-pres}, the structures $\p{M,0,1,+,<}$ and $\p{P_2,1,2,\cdot,<}$ are both
models of Presburger arithmetic, hence elementarily equivalent, hence isomorphic by Lemma~\ref{lem:rec-sat}. If $2^x$ is
such an isomorphism, then $\p{\sM,2^x}\model\teipe$.
\end{Pf}

The growth axiom interconnects the structures $\p{M,0,1,+,<}$ and $\p{P_2,1,2,\cdot,<}$ in a way that seems to preclude
a similarly easy proof of the extension of Theorem~\ref{thm:teipp} to $\teipe^+$. One idea that does not work is to use the
joint consistency theorem to expand an elementary extension of~$\sM$ to a model of $\Th(\p{\N,0,1,+,<,2^x})$, taking
$\Th(\p{\N,0,1,+,<,\cdot\res P_2})$ as the common subtheory: since $\Th(\p{\N,0,1,+,<,2^x})$ is decidable due to
Semenov~\cite{sem:plus-pow}, its reduct $\Th(\p{\N,0,1,+,<,\cdot\res P_2})$ is recursively axiomatizable, thus in
principle it might be possible to take its axiomatization by a few natural axioms or schemata and check that it is
included in~$\teipp$. This fails for two reasons: first, even though $\Th(\p{\N,0,1,+,<,2^x})$ is explicitly
axiomatized in Cherlin and Point~\cite{cher-poi:presext}, we do not know of a similar axiomatization of
$\Th(\p{\N,0,1,+,<,\cdot\res P_2})$ in the literature, and this would likely require some work to devise. Second,
$\teipe^+$ does \emph{not}, in fact, include even the weaker theory $\Th(\p{\N,0,1,+,<,P_2})$: as we will see below, it
does not prove that powers of~$2$ are not divisible by~$3$.

In absence of a better idea, we will get our hands dirty and construct the required $2^x$ obeying GA by
a back-and-forth argument (cf.\ \cite[Thm.~6.4]{ej:vtceip}):
\begin{Thm}\label{thm:teipp-ga}
The first-order theory of EIP of RCEF satisfying GA in $\langp$ is $\teipp$. 
\end{Thm}
\begin{Pf}
It suffices to show that $\teipe^+$ is conservative over $\teipp$. Let $\sM=\p{M,0,1,+,\cdot,{<},P_2}$ be a countable
recursively saturated model of $\teipp$; we will show that it expands to a model of~$\teipe^+$. Let $\sR$ be a RCF
whose IP is $\sM_\pm$, $P_2^{\pm1}=\{u,u^{-1}:u\in P_2\}\sset R_{>0}$, $B=\{x\in R:\exists n\in\N\,\abs x\le n\}$, and
$B^\times=\{x\in R:\exists{n\in\N}\,{n^{-1}\le x\le n}\}$. If $u,v\in P_2^{\pm1}$ and $m\in\N_{>0}$,
$u\equiv^\times v\pmod m$ means $uv^{-1}=w^m$ for some $w\in P_2^{\pm1}$; by Lemma~\ref{lem:teipp-pres},
$u\equiv^\times2^l\pmod m$ for a unique $l$ such that $0\le l<m$.

Fix enumerations $M=\p{c_i:i\in\omega}$ and $P_2=\p{d_i:i\in\omega}$. We will construct sequences
$\p{a_i:i\in\omega}\sset M$ and $\p{b_i:i\in\omega}\sset P_2$ so that they satisfy the following properties for each
$k\ge1$ by induction on~$k$:
\begin{enumerate}
\item\label{item:1} $a_0=1$, $b_0=2$, $a_{2i+2}=c_i$, $b_{2i+1}=d_i$.
\item\label{item:2} For all $\vec q\in\Q^k$, $\sum_{i<k}q_ia_i=0\implies\prod_{i<k}b_i^{q_i}=1$. (Here,
$b_i^{q_i}\in R_{>0}$.)
\item\label{item:3} For all $0\le l<m\in\N$ and $i<k$, $a_i\equiv l\pmod m\implies b_i\equiv^\times2^l\pmod m$.
\item\label{item:4} For all $\vec q\in\Q^k$, $\prod_{i<k}b_i^{q_i}>\sum_{i<k}q_ia_i$.
\end{enumerate}
We observe that conditions \ref{item:2} and~\ref{item:3} can be stated with $\iff$ in place of $\implies$.
For~\ref{item:3}, this follows from the uniqueness of $l<m$ such that $a_i\equiv l\pmod m$, resp.\
$b_i\equiv^\times2^l\pmod m$. For~\ref{item:2}, this follows from \ref{item:4}: if $\sum_iq_ia_i\ne0$, then
$\sum_inq_ia_i\ge1$ for some $n\in\Z$, thus $\bigl(\prod_ib_i^{q_i}\bigr)^n>1$ by~\ref{item:4}, and in particular,
$\prod_ib_i^{q_i}\ne1$. The same argument actually shows that the conditions imply that the
mapping $a_i\mapsto b_i$ extends to an ordered $\Q$-linear space isomorphism between the span of $\{a_i:i<k\}$ in $\p{R,+,<}$, and the span of $\{b_i:i<k\}$ in $\p{R_{>0},\cdot,<}$, i.e.,
\begin{equation}\label{eq:17}
\sum_{i<k}q_ia_i>0\iff\prod_{i<k}b_i^{q_i}>1;
\end{equation}
likewise,
\begin{equation}\label{eq:19}
\sum_{i<k}q_ia_i>\N\iff\prod_{i<k}b_i^{q_i}>\N\colon
\end{equation}
the left-to-right implication follows from \ref{item:4}, while if $\sum_{i<k}q_ia_i<n$ for some~$n\in\N$, then
\ref{item:4} applied to $na_0-\sum_iq_ia_i>0$ gives $\prod_ib_i^{q_i}<2^n$.

We also observe that condition~\ref{item:4} is equivalent to
\begin{equation}\label{eq:18}
\sum_{i<k}q_ia_i>\N\implies\prod_{i<k}b_i^{q_i}>\sum_{i<k}q_ia_i
\end{equation}
for all $\vec q\in\Q^k$, as other conditions imply the conclusion when $\sum_iq_ia_i>\N$ does not hold: let
$r=\sum_iq_ia_i$. If $r\le0$, there is nothing to prove, as ${\prod_{i<k}b_i^{q_i}>0}$. If $0\le r\in B$, then $r\in\Q$
(if $q_i=n_i/m$ for some $\vec n,m\in\Z$, $m>0$, then $mr=\sum_in_ia_i\in M\cap B=\Z$ as $\sM$ is a model of Presburger
arithmetic). Then in view of~\ref{item:1}, $\sum_iq_ia_i-ra_0=0$ implies $b_0^{-r}\prod_ib_i^{q_i}=1$ by~\ref{item:2},
i.e., $\prod_ib_i^{q_i}=2^r>r$ (referring to the standard exponential).

It is clear that after we finish the construction, the conditions ensure that $a_i\mapsto b_i$ defines an isomorphism
$2^x\colon\p{M,0,1,+,<}\to\p{P_2,1,2,\cdot,<}$, and $\p{\sM,2^x}\models\teipe^+$.

We now proceed with the construction. For $k=1$, we put $a_0=1$, $b_0=2$ as requested by~\ref{item:1}; then
\ref{item:2}--\ref{item:4} hold. Having constructed $\p{a_i:i<k}$ and~$\p{b_i:i<k}$ satisfying
\ref{item:2}--\ref{item:4}, we will construct $a_k$ and~$b_k$ as follows.

Assume that $k$ is even. Put $a_k=c_{k/2-1}$; we need to find a matching~$b_k\in P_2$. First, if
$a_k+\sum_{i<k}q_ia_i\in B$ for some $\vec q\in\Q^k$, then $a_k=\sum_{i<k}q_ia_i$ for some $\vec q\in\Q^k$ by the same
argument as in the equivalence of \ref{item:4} and~\eqref{eq:18} above, and we define $b_k=\prod_{i<k}b_i^{q_i}$. Write
$q_i=n_i/m$ for some $\vec n\in\Z^k$ and $m\in\N_{>0}$, and let $0\le l_i<m$ be such that $a_i\equiv l_i\pmod m$. Then
$0\equiv ma_k\equiv\sum_in_il_i\pmod m$. Using \ref{item:3} from the induction hypothesis,
$b_i\equiv^\times2^{l_i}\pmod m$, thus $\prod_ib_i^{n_i}\equiv^\times2^{\sum_in_il_i}\equiv^\times1\pmod m$. This shows
that $b_k=\bigl(\prod_ib_i^{n_i}\bigr)^{1/m}\in P_2$; moreover, an analogous argument gives~\ref{item:3}. Conditions
\ref{item:2} and~\ref{item:4} follow from the induction hypothesis.

Now, assume that $a_k+\sum_iq_ia_i\notin B$ for all $\vec q\in\Q^k$. Then condition \ref{item:2} will follow from the
induction hypothesis for whatever choice of~$b_k$, hence we only need $b_k$ to satisfy \ref{item:3} and~\eqref{eq:18}.
Condition~\eqref{eq:18} for $q_k=0$ follows from the induction hypothesis. For $q_k>0$, the condition
\begin{equation}\label{eq:14}
q_ka_k+\sum_{i<k}q_ia_i>\N\implies b_k^{q_k}\prod_{i<k}b_i^{q_i}>q_ka_k+\sum_{i<k}q_ia_i
\end{equation}
is equivalent to
\begin{equation}\label{eq:12}
a_k>\sum_{i<k}r_ia_i\implies b_k\prod_{i<k}b_i^{-r_i}>\Bigl[q_k\Bigl(a_k-\sum_{i<k}r_ia_i\Bigr)\Bigr]^{1/q_k},
\end{equation}
where $r_i=-q_i/q_k$ (using that $a_k-\sum_{i<k}r_ia_i>0$ implies $a_k-\sum_{i<k}r_ia_i>\N$). Also, we have
$q_k\Bigl(a_k-\sum_ir_ia_i\Bigr)<\Bigl(a_k-\sum_ir_ia_i\Bigr)^2$, thus \eqref{eq:12} holds for all $q_k>0$
and $\vec q\in\Q^k$ iff
\begin{equation}\label{eq:13}
a_k>\sum_{i<k}r_ia_i\implies b_k>\Bigl(a_k-\sum_{i<k}r_ia_i\Bigr)^n\prod_{i<k}b_i^{r_i}
\end{equation}
holds for all $\vec r\in\Q^k$ and $n\in\N$. Likewise, \eqref{eq:14} for all $q_k<0$ and $\vec q\in\Q^k$ is equivalent
to
\begin{equation}\label{eq:15}
a_k<\sum_{i<k}q_ia_i\implies b_k<\Bigl(\sum_{i<k}q_ia_i-a_k\Bigr)^{-n}\prod_{i<k}b_i^{q_i}
\end{equation}
for all $\vec q\in\Q^k$ and $n\in\N$. Thus, to satisfy conditions \ref{item:2}--\ref{item:4}, it is enough to take for
$b_k$ any realizer of the type
\begin{align*}
\Gamma(x)=\{P_2(x)\}&\cup\bigl\{a_k\equiv l\pmod m\to x\equiv^\times2^l\pmod m:0\le l<m\in\N\bigr\}\\
&\cup\Bigl\{a_k>\sum_{i<k}r_ia_i\to x>\Bigl(a_k-\sum_{i<k}r_ia_i\Bigr)^n\prod_{i<k}b_i^{r_i}:\vec r\in\Q^k,n\in\N\Bigr\}\\
&\cup\Bigl\{a_k<\sum_{i<k}q_ia_i\to x<\Bigl(\sum_{i<k}q_ia_i-a_k\Bigr)^{-n}\prod_{i<k}b_i^{q_i}:\vec q\in\Q^k,n\in\N\Bigr\}.
\end{align*}
Observe that $\Gamma(x)$ can indeed be expressed as a recursive type in $\langp$ with parameters $\vec a,\vec b$:
e.g., if $r_i=n_i/m$ with $\vec n\in\Z^k$ and $m\in\N_{>0}$, then $x>(\cdots)^n\prod_ib_i^{r_i}$ is equivalent to
$x^m\prod_{n_i<0}b_i^{-n_i}>(\cdots)^{nm}\prod_{n_i>0}b_i^{n_i}$, etc. Thus, using the recursive saturation of~$\sM$, it
only remains to check that every finite $\Gamma_0\sset\Gamma$ is satisfiable.

Apart from $P_2(x)$, the formulas in $\Gamma_0$ are implications whose premises do not depend on~$x$; we may discard
those whose premises are false, and simplify the remaining ones by removing their premises. If $a_k\equiv l_j\pmod{m_j}$
for $j<t$, then $a_k\equiv l\pmod m$, where $m=\lcm(m_0,\dots,m_{t-1})$ and $l\equiv l_j\pmod{m_j}$; then
$x\equiv^\times2^l\pmod m$ implies $x\equiv^\times2^{l_j}\pmod{m_j}$ for each $j<t$. Thus, we may assume $\Gamma_0$
contains only one congruence $x\equiv^\times2^l\pmod m$. Likewise, we can take the maximum (minimum) right-hand side
among the inequalities $x>\cdots$ ($x<\cdots$, resp.), thus we may assume that $\Gamma_0$ contains one inequality
of the form $x>\cdots$ (we may assume there is at least one by considering e.g.\ $\vec r=\vec0$ and $n=0$, which gives
$x>1$), and at most one inequality of the form $x<\cdots$. If there is no inequality $x<\cdots$, it is easy to see that
$\Gamma_0$ is satisfiable, hence we may assume that
\[\Gamma_0=\Bigl\{P_2(x),x\equiv^\times2^l\pmod m,\Bigl(a_k-\sum_{i<k}r_ia_i\Bigr)^n\prod_{i<k}b_i^{r_i}<x,
  x<\Bigl(\sum_{i<k}q_ia_i-a_k\Bigr)^{-n}\prod_{i<k}b_i^{q_i}\Bigr\}\]
for some $0\le l<m\in\N$, $\vec q,\vec r\in\Q^k$, and $n\in\N$, where
\[\sum_{i<k}r_ia_i<a_k<\sum_{i<k}q_ia_i.\]
(We may assume both inequalities use the same~$n$ by enlarging one if necessary.) Since $P_2^\pm$ is an IP of
$\p{R_{>0},1,\cdot,<}$ (cf.\ Lemma~\ref{lem:teipp-pres}), there exists an element $x\in P_2$ satisfying
$x\equiv^\times2^l\pmod m$ in any interval $[u,v)$ such that $v\ge2^mu>0$. Thus, $\Gamma_0$ is satisfiable if
\[\Bigl(\sum_{i<k}q_ia_i-a_k\Bigr)^{-n}\prod_{i<k}b_i^{q_i}>2^m\Bigl(a_k-\sum_{i<k}r_ia_i\Bigr)^n\prod_{i<k}b_i^{r_i},\]
i.e.,
\begin{equation}\label{eq:16}
\prod_{i<k}b_i^{q_i-r_i}>2^m\Bigl(a_k-\sum_{i<k}r_ia_i\Bigr)^n\Bigl(\sum_{i<k}q_ia_i-a_k\Bigr)^n.
\end{equation}
Now, using $\sum_i(q_i-r_i)a_i>\N$, we have
\[\Bigl(\frac1{2n+1}\sum_{i<k}(q_i-r_i)a_i\Bigr)^{2n+1}
  >2^m\Bigl(a_k-\sum_{i<k}r_ia_i\Bigr)^n\Bigl(\sum_{i<k}q_ia_i-a_k\Bigr)^n,\]
whence \eqref{eq:16} follows from the instance
\[\prod_{i<k}b_i^{(q_i-r_i)/(2n+1)}>\sum_{i<k}\frac{q_i-r_i}{2n+1}a_i\]
of the induction hypothesis. This finishes the construction of $a_k$ and~$b_k$ for $k$ even.

Let $k$ be odd, and put $b_k=d_{(k-1)/2}$; we will find a suitable~$a_k$. If $b_k\prod_{i<k}b_i^{q_i}\in B^\times$ for
some $\vec q\in\Q^k$, then as in the case of even~$k$, we obtain $b_k=\prod_{i<k}b_i^{q_i}$ for some $\vec q\in\Q^k$, and
then $a_k=\sum_{i<k}q_ia_i$ will satisfy \ref{item:2}--\ref{item:4}. Thus, we may assume
\begin{equation}\label{eq:24}
b_k\prod_{i<k}b_i^{q_i}\notin B^\times
\end{equation}
for all $\vec q\in\Q^k$. Then \ref{item:2} and~\ref{item:4} will hold if $a_k$ satisfies
\[b_k^{q_k}\prod_{i<k}b_i^{q_i}>1\implies b_k^{q_k}\prod_{i<k}b_i^{q_i}>q_ka_k+\sum_{i<k}q_ia_i>0\]
for all $q_k\ne0$ and $\vec q\in\Q^k$. Similarly to the case of even~$k$, one can check that this amounts to the
conditions
\begin{align*}
b_k>\prod_{i<k}b_i^{q_i}&\implies\sum_{i<k}q_ia_i<a_k<\sum_{i<k}q_ia_i+\Bigl(b_k\prod_{i<k}b_i^{-q_i}\Bigr)^{1/n},\\
b_k<\prod_{i<k}b_i^{q_i}&\implies\sum_{i<k}q_ia_i-\Bigl(b_k^{-1}\prod_{i<k}b_i^{q_i}\Bigr)^{1/n}<a_k<\sum_{i<k}q_ia_i
\end{align*}
for all $\vec q\in\Q^k$ and $n\in\N_{>0}$. Thus, using recursive saturation, it suffices to show that each finite
subset $\Gamma_0$ of the type
\begin{align*}
\Gamma(x)&=\bigl\{x\equiv l\pmod m:0\le l<m\in\N,b_k\equiv^\times2^l\pmod m\bigr\}\\
&\qquad\cup\Bigl\{x>\sum_{i<k}q_ia_i:\vec q\in\Q^k,b_k>\prod_{i<k}b_i^{q_i}\Bigr\}\\
&\qquad\cup\Bigl\{x<\sum_{i<k}r_ia_i:\vec r\in\Q^k,b_k<\prod_{i<k}b_i^{r_i}\Bigr\}\\
&\qquad\cup\Bigl\{x>\sum_{i<k}s_ia_i-\Bigl(b_k^{-1}\prod_{i<k}b_i^{s_i}\Bigr)^{1/n}:\vec s\in\Q^k,n\in\N_{>0},b_k<\prod_{i<k}b_i^{s_i}\Bigr\}\\
&\qquad\cup\Bigl\{x<\sum_{i<k}t_ia_i+\Bigl(b_k\prod_{i<k}b_i^{-t_i}\Bigr)^{1/n}:\vec t\in\Q^k,n\in\N_{>0},b_k>\prod_{i<k}b_i^{t_i}\Bigr\}
\end{align*}
is satisfiable. (To make the type recursive, we would write it with implications as in the case of even~$k$.) Again, we
may assume that $\Gamma_0$ consists of one congruence $x\equiv l\pmod m$, one lower bound on~$x$, and one upper bound;
it will be satisfiable as long as the difference between the upper and lower bounds is larger than~$m$.
Thus, assume that
\[\prod_{i<k}b_i^{q_i},\prod_{i<k}b_i^{t_i}<b_k<\prod_{i<k}b_i^{r_i},\prod_{i<k}b_i^{s_i};\]
we need to check that
\begin{align}
\label{eq:20}\sum_{i<k}q_ia_i+m&<\sum_{i<k}r_ia_i,\\
\label{eq:21}\sum_{i<k}q_ia_i+m&<\sum_{i<k}t_ia_i+\Bigl(b_k\prod_{i<k}b_i^{-t_i}\Bigr)^{1/n},\\
\label{eq:22}\sum_{i<k}s_ia_i-\Bigl(b_k^{-1}\prod_{i<k}b_i^{s_i}\Bigr)^{1/n}+m&<\sum_{i<k}r_ia_i,\\
\label{eq:23}\sum_{i<k}s_ia_i-\Bigl(b_k^{-1}\prod_{i<k}b_i^{s_i}\Bigr)^{1/n}+m&<\sum_{i<k}t_ia_i+\Bigl(b_k\prod_{i<k}b_i^{-t_i}\Bigr)^{1/n}.
\end{align}

Since $\prod_ib_i^{r_i-q_i}>\N$ by~\eqref{eq:24}, the inequality~\eqref{eq:20} follows from~\eqref{eq:19}. For~\eqref{eq:21}, we have
\[\sum_{i<k}(q_i-t_i)a_i+m=n\Bigl(\sum_{i<k}\frac{q_i-t_i}na_i+\frac mna_0\Bigr)
 <n2^{m/n}\prod_{i<k}b_i^{(q_i-t_i)/n}<\Bigl(b_k\prod_{i<k}b_i^{-t_i}\Bigr)^{1/n},\]
using \ref{item:4} and $b_k>n^n2^m\prod_ib_i^{q_i}$ from~\eqref{eq:24}; the argument for~\eqref{eq:22} is similar.
Finally,
\begin{align*}
\sum_{i<k}(s_i-t_i)a_i+m&=(2n+1)\Bigl(\sum_{i<k}\frac{s_i-t_i}{2n+1}a_i+\frac m{2n+1}\Bigr)\\
&<(2n+1)2^{m/(2n+1)}\prod_{i<k}b_i^{(s_i-t_i)/(2n+1)}\\
&<\prod_{i<k}b_i^{(s_i-t_i)/(2n)}\le\Bigl(b_k^{-1}\prod_{i<k}b_i^{s_i}\Bigr)^{1/n}+\Bigl(b_k\prod_{i<k}b_i^{-t_i}\Bigr)^{1/n}
\end{align*}
using \ref{item:4},
\[\prod_{i<k}b_i^{s_i-t_i}=\Bigl(b_k^{-1}\prod_{i<k}b_i^{s_i}\Bigr)\Bigl(b_k\prod_{i<k}b_i^{-t_i}\Bigr)
\le\max\Bigl(\Bigl\{b_k^{-1}\prod_{i<k}b_i^{s_i},b_k\prod_{i<k}b_i^{-t_i}\Bigr\}\Bigr)^2,\]
and $\prod_{i<k}b_i^{s_i-t_i}>\N$, which follows from \eqref{eq:24}.
\end{Pf}

\begin{Exm}\label{exm:shep}
There exists a countable model of $\io$ that expands to a model of $\teipp$, but not to a model of $\teipe$.
\end{Exm}
\begin{Pf}
Let $\sM_\pm$ be the ring of Puiseux polynomials $\sum_{q\in Q}a_qx^q$ with $Q\sset\Q_{\ge0}$ finite, $a_q$ real
algebraic, and $a_0\in\Z$, ordered so that $x>\N$. Its nonnegative part $\sM$ is a model of $\io$ by
Shepherdson~\cite{sheph}, and it can be checked readily that $\p{\sM,P_2}\model\teipp$, where
\[P_2=\{2^nx^q:q\in\Q_{\ge0},n\in\Z,(q>0\text{ or }n\ge0)\}.\] On the other hand, assume for contradiction that
$\p{\sM,2^x}\model\teipe$. Then $2^x$ extends to an ordered group embedding
$2^x\colon\sM_\pm\to\sF^\times_{>0}$, where $\sF^\times_{>0}=\p{F_{>0},1,\cdot,<}$ is the multiplicative
group of positive elements of the fraction field $\sF$ of~$\sM_\pm$. Since the image of $2^x$ is an IP of
$\sF^\times_{>0}$, $2^x$ induces an isomorphism of the ordered groups $\sM_\pm/\Z$ and $\sF^\times_{>0}/B^\times$,
where $B^\times=\{x\in F_{>0}:\exists{n\in\N}\,{n^{-1}\le x\le n}\}$. But every coset of $B^\times$ contains exactly
one monomial $x^q$, $q\in\Q$, thus $\sF^\times_{>0}/B^\times\simeq\p{\Q,0,+,<}$ is archimedean, whereas $\sM_\pm/\Z$,
isomorphic to the additive group of Puiseux polynomials with $a_0=0$, is nonarchimedean. This is a contradiction.
\end{Pf}

\section{Exponential integer parts in $\langor$}\label{sec:teip}

We now turn to the most interesting case, namely the theory of EIP of RCEF in the basic language of
arithmetic~$\langor$. Our axiomatization of this theory will express the existence of winning strategies in a certain
game on integers. We describe the game first to motivate the definition of the theory.
\begin{Def}\label{def:powg}
Let $\sM\model\io$ and $\alpha\le\omega$. The \emph{power-of-$2$ game $\powg_\alpha(\sM)$} is played between two
players, \emph{Challenger} (C) and \emph{Powerator} (P), in $\alpha$~rounds: in each round $0\le i<\alpha$, C picks
$x_i\in M_{>0}$, and P responds with $u_i\in M_{>0}$ such that $u_i\le x_i<2u_i$. C wins the game if
$u_iu_j<u_h<2u_iu_j$ for some $h,i,j<\alpha$, otherwise P wins.

More generally, if $t\le\alpha$ is finite, and $u_0,\dots,u_{t-1}\in M_{>0}$, let $\powg^t_\alpha(\sM,\vec u)$ denote
the $\powg_\alpha(\sM)$ game where the first $t$ responses by P are fixed as $\vec u$ (the values of $x_i$, $i<t$, do
not matter, as they do not enter the winning condition; for definiteness, we may imagine $x_i=u_i$). We may write just
$\powg^t_\alpha(\vec u)$ if $\sM$ is understood from the context.
\end{Def}

While not being part of the official rules as we want to keep them simple, we will often use the following
alternative conditions:
\begin{Obs}\label{obs:powg-div}
\ 
\begin{enumerate}
\item\label{item:7} If $u_i\le u_j$ and $u_i\nmid u_j$ for some $i,j<h<\alpha$, then Challenger can win the game in
round $h$ by playing $x_h=\fl{u_j/u_i}$.
\item\label{item:8} For any $x_i>0$, Challenger can force Powerator to respond with $u_i$ such that $x_i\le u_i<2x_i$.
\end{enumerate}
\end{Obs}
\begin{Pf}
\ref{item:7}: P must respond with $u_h$ such that $u_h\le x_h<u_j/u_i<2u_h$, i.e., $u_hu_i<u_j<2u_hu_i$.

\ref{item:8}: Let C play $2x_i-1$, so that $u_i\le2x_i-1<2u_i$.
\end{Pf}

\begin{Rem}\label{rem:powg-even}
C cannot go wrong by restricting their moves to even numbers: instead of playing $2x+1$, to which the valid responses
of P are in $\{x+1,\dots,2x+1\}$, C can play $2x$ with valid responses in $\{x+1,\dots,2x\}$, unless $x=0$. A move
$x_i=1$, forcing P to reply with $u_i=1$, can be eliminated as well: let C skip the move. The only way this can affect
the game is when we reach a position with $u_j<u_h<2u_j$ for some $h,j$ (which would make C win as
$1\cdot u_j<u_h<2\cdot1\cdot u_j$); then C can play $x_l=2(u_ju_h-1)$ en lieu of the skipped round, forcing P to reply
with $u_ju_h\le u_l<2u_ju_h$ and lose, as either $u_ju_h<u_l<2u_ju_h$ or $u_j^2<u_l=u_ju_h<2u_j^2$.
\end{Rem}

The intuition behind the game is that Powerator can win by playing powers of~$2$:
\begin{Lem}\label{lem:powg-teipp}
If $\p{\sM,P_2}\model\teipp$, then Powerator has a winning strategy in $\powg_\alpha(\sM)$ for every $\alpha\le\omega$,
and more generally, in $\powg^t_\alpha(\sM,\vec u)$ for every $t<\omega$, $t\le\alpha$, and $\vec u\sset P_2$.
\end{Lem}
\begin{Pf}
By Lemmas \ref{lem:teipp} and~\ref{lem:teippp}, $\p{\sM,P_2}\model\teipp''$. Given a move $x_i$ of C, let P respond with
$u_i\in P_2$ such that $u_i\le x_i<2u_i$, which exists by~\eqref{eq:4}. Then $u_iu_j<u_h<2u_iu_j$ is impossible
by~\eqref{eq:11}.
\end{Pf}

\begin{Def}\label{def:teip}
For any $t\le n<\omega$, let $\pwin^t_n(u_0,\dots,u_{t-1})$ denote the formula
\[\forall x_t\:\exists u_t\dots\forall x_{n-1}\:\exists u_{n-1}
\Bigl(\ET_{t\le i<n}(x_i>0\to u_i\le x_i<2u_i)\land\ET_{\mclap{h,i,j<n}}\neg(u_iu_j<u_h<2u_iu_j)\Bigr),\]
expressing that Powerator has a winning strategy in $\powg^t_n(\sM,\vec u)$.

$\teip$ is the $\langor$-theory axiomatized by $\io+\{\pwin^0_n:n\in\omega\}$.
\end{Def}

The basic properties below follow immediately from the definition:
\begin{Lem}\label{lem:theta}
If $t<n$, then $\pwin^t_n(\vec u)$ is equivalent to
\[\forall x_t>0\,\exists u_t\,\bigl(u_t\le x_t<2u_t\land\pwin^{t+1}_n(\vec u,u_t)\bigr).\]
If $t\le m<n$, then $\pwin^t_n(\vec u)$ implies $\pwin^t_m(\vec u)$.
\noproof\end{Lem}
\pagebreak[2]
\begin{Thm}\label{thm:teip}
The first-order theory of EIP of RCEF in $\langor$, with or without GA, is $\teip$.
\end{Thm}
\begin{Pf}
In view of Theorems \ref{thm:teipp} and~\ref{thm:teipp-ga}, it suffices to show that $\teipp$ is a conservative extension of~$\teip$.
Clearly, $\sM\model\io$ is a model of $\teip$ iff Powerator has a winning strategy in $\powg_n(\sM)$ for all
$n\in\omega$; in particular, $\teipp\vdash\teip$ follows from Lemma~\ref{lem:powg-teipp}.

On the other hand, let $\sM$ be a countable recursively saturated model of~$\teip$; we will expand $\sM$ to a model of
$\teipp$. The basic idea is that due to recursive saturation, P also has a winning strategy in $\powg_\omega(\sM)$, and then
if we let C enumerate $M_{>0}$, the responses of P form a set $P_2$ such that $\p{\sM,P_2}\model\teipp''$.

Formally, let $\Gamma^t(u_0,\dots,u_{t-1})=\{\pwin^t_n(\vec u):t\le n<\omega\}$ for $t<\omega$, and fix an enumeration
$\p{a_i:i<\omega}$ of~$M_{>0}$. We will construct a sequence $\p{b_i:i<\omega}\sset M_{>0}$ such that
$b_i\le a_i<2b_i$ and $\sM\model\Gamma^t(\vec b)$ by induction on~$t$.

We have $\sM\model\Gamma^0$ as $\Gamma^0\sset\teip$. Assuming $\sM\model\Gamma^t(b_0,\dots,b_{t-1})$, we can take for
$b_t$ any realizer of the type $\Gamma^{t+1}(\vec b,u_t)\cup\{u_t\le a_t<2u_t\}$, hence using recursive saturation, we
only need to check its finite satisfiability. In view of Lemma~\ref{lem:theta}, it suffices to observe that
for any $n>t$, $\sM\model\exists u_t\,\bigl(u_t\le a_t<2u_t\land\pwin^{t+1}_n(\vec b,u_t)\bigr)$ follows from
$\sM\model\pwin^t_n(\vec b)$.

When the construction of $\p{b_i:i<\omega}$ is finished, let $P_2=\{b_i:i<\omega\}$. Then the properties of $\vec a$
and $\vec b$ ensure $\p{\sM,P_2}\model\teipp''$.
\end{Pf}

Coupled with Lemma~\ref{lem:powg-teipp}, the proof gives a characterization of $\langor$-reducts of countable
models of $\teipp$:
\begin{Cor}\label{cor:ctbl-powg-w}
A countable model $\sM\model\io$ expands to a model of $\teipp$ iff Powerator has a winning strategy in
$\powg_\omega(\sM)$.
\noproof\end{Cor}

Using our axiomatization of $\teip$, it is now easy to answer negatively Question~7.3 from~\cite{ej:vtceip}.
\begin{Exm}\label{exm:teip-io}
The following consequence of $\teip$ is not provable in $\io$:
\begin{equation}\label{eq:26}
\forall x\:\exists u\ge x\:\forall y\:\bigl(0<y<x\to\exists v\:(v\le y<2v\land v\mid u)\bigr).
\end{equation}
(We can make it $\Pi_1$ by further bounding $u<2x$.) Thus, some models of $\io$ have no elementary extension to an EIP
of a RCEF.
\end{Exm}
\begin{Pf}
First, \eqref{eq:26} indeed follows from $\teip$ (specifically, $\pwin^0_3$) in view of Observation~\ref{obs:powg-div}.

On the other hand, Smith~\cite{smith} constructed a nonstandard $\sM\model\io$ which is a UFD (or even PID): i.e.,
every $x\in M_{>0}$ can be written as a product of a sequence $\p{p_i:i<k}$ of primes $p_i\in M_{>0}$ of
standard length $k\in\N$. It follows that $x^*=\prod_{i:p_i\in\N}p_i$ is the largest standard divisor of~$x$.

Assume for contradiction that \eqref{eq:26} holds in~$\sM$. Let $x\in M$ be nonstandard, and $u\in M$ satisfy the
conclusion of \eqref{eq:26}. Take $y=2u^*$ (which is standard, thus $y<x$), and let $v\in M$ be such that $v\le y<2v$
and $v\mid u$. Then $v$ is a standard divisor of~$u$, but $v>u^*$, a contradiction.
\end{Pf}

We mention that $\io\vdash\pwin^0_2$, hence the use of $\pwin^0_3$ in Example~\ref{exm:teip-io} is the best possible.

\section{Analysis of $\powg$}\label{sec:analysis-powg}

Unlike $\teipe$ and $\teipp$, we defined $\teip$ by an infinite axiom schema, but it is not clear whether this is
necessary:
\begin{Que}\label{que:finax}
Is $\teip$ finitely axiomatizable over $\io$?
\end{Que}

We do not know how to resolve this question, but we can at least give a partial answer. Let us observe that if there
were only finitely many inequivalent formulas among $\{\pwin^1_n:n\in\omega\}$, then $\teip$ would be finitely
axiomatizable over~$\io$ by Lemma~\ref{lem:theta}. We will show that this is not the case, though:
the formulas $\{\pwin^1_n:n\in\omega\}$ are strictly increasing in strength, even over $\Th(\N)$. This is equivalent to
$\{c(u):u\in\N_{>0}\bez P_2^\N\}=\N_{>0}$, using the notation below:
\begin{Def}\label{def:cxity}
If $\sM\model\io$ and $\vec u\in M_{>0}^t$, the \emph{$\powg$-complexity} of $\vec u$, denoted $c(\sM,\vec u)$, is the
least $n\in\omega$ such that C has a winning strategy in $\powg^t_{t+n}(\sM,\vec u)$; if such an $n$ does not exist, we
put $c(\sM,\vec u)=\infty$. If $\sM=\N$, we write just $c(\vec u)$. (Observe that $c(u)\ge n$ iff
$\N\model\pwin^1_n(u)$.) Let $P_2^\N$ denote the set of powers of~$2$
in~$\N$.
\end{Def}
\pagebreak[2]
\begin{Lem}\label{lem:c-def}
For any $\vec u\in\N_{>0}^t$, $c(\vec u)$ is finite iff some $u_i$ is not a power of~$2$.
\end{Lem}
\begin{Pf}
The left-to-right implication follows from Lemma~\ref{lem:powg-teipp}. On the other hand, if $c(\vec u)=\infty$,
i.e., P has a winning strategy in $\powg^t_n(\vec u)$ for all $n\ge t$, and $\sM$ is a countable recursively saturated
model of~$\Th(\N)$, then $\sM$ expands to a model $\p{\sM,P_2}\model\teipp$ such that $\vec u\sset P_2$ by
the proof of Theorem~\ref{thm:teip}. Taking one more elementary extension if necessary, it expands to a model
$\p{\sM,2^x}\model\teipe$ such that $\vec u\sset\im(2^x)$. But $\teipe$ implies that $2^x$ extends the standard
function and maps nonstandard values to nonstandard values, hence $\vec u\sset P_2^\N$.

(The reader is invited to construct a simple explicit winning strategy for C if some $u_i$ is not a power of~$2$.
We will present an optimized one below in Theorem~\ref{thm:powg-ub}.)
\end{Pf}

For the application to finite non-axiomatizability of $\{\pwin^1_n(u):n\ge1\}$, it would be clearly enough to show that
$\{c(u):u\in\N_{>0}\bez P_2^\N\}$ is unbounded. Let us observe that this is, in fact, equivalent to
$\{c(u):u\in\N_{>0}\bez P_2^\N\}=\N_{>0}$:
\begin{Lem}\label{lem:c-init}
$\{c(u):u\in\N_{>0}\bez P_2^\N\}$ is an initial segment of $\N_{>0}$.
\end{Lem}
\begin{Pf}
There are $u$ such that $c(u)=1$, see Example~\ref{exm:c1}. Let $n>1$, and assume that there exists $u\notin P_2^\N$ such that $c(u)\ge
n$. Let $u$ be the smallest such number; we will give a strategy for C showing $c(u)=n$. First, C plays $u-1$, thus P
responds with a $v$ such that $u/2\le v<u$. If $u/2<v$, then $v\nmid u$, thus C can win in the second round by
Observation~\ref{obs:powg-div}; otherwise, $v=u/2\notin P_2^\N$, thus $c(u/2)<n$ by the minimality of~$u$, and C can just follow
the optimal strategy for $u/2$.
\end{Pf}

While our goal is to prove lower bounds on $c(u)$, we will start with upper bounds to get an idea of what is in the
realm of possible:
it turns out that C can very efficiently exploit irregularities in exponents of the prime factorization of~$u$, hence
our lower bounds will need to be somewhat delicate. The main tool of Challenger is the following divide-and-conquer
strategy.
\begin{Lem}\label{lem:div-conq}
Let $n>1$ and $u,v,\vec u\in\N_{>0}$.
\begin{enumerate}
\item\label{item:9} If $u^n<v<2u^n$, then $c(u,v)\le\cl{\log\fl{\log n}}+1$.
\item\label{item:5} $c(\vec u,u)\le\max\bigl\{c(\vec u,u,u^n)+1,\cl{\log\fl{\log n}}+2\bigr\}$.
\end{enumerate}
\end{Lem}
\begin{Pf}
\ref{item:9}: Let $k=\fl{\log n}$, and put $i_0=0$, $j_0=k$, $v_0=v$ so that $u^{\tdive n{i_0}}<v_0<2u^{\tdive n{i_0}}$
and $u^{\tdive n{j_0}}=u$. Using Observation~\ref{obs:powg-div}, let C play $2u^{\tdive n{m_1}}-1$ for
$m_1=\fl{k/2}=\fl{(i_0+j_0)/2}$ so that P responds with a $w_1$ such that $u^{\tdive n{m_1}}\le w_1<2u^{\tdive
n{m_1}}$. If $w_1=u^{\tdive n{m_1}}$, put $i_1=i_0$, $j_1=m_1$, and $v_1=v_0$; otherwise, $i_1=m_1$, $j_1=j_0$, and
$v_1=w_1$. Either way, the responses of P include $u^{\tdive n{j_1}}$ and $v_1$ satisfying $u^{\tdive
n{i_1}}<v_1<2u^{\tdive n{i_1}}$, where $i_1<j_1$ and $j_1-i_1\le\cl{k/2}$. We continue a binary search in the same
way: after $\cl{\log k}$ rounds, the responses of P will include $u'=u^{\tdive n{i+1}}$ and $v'$ such that $u^{\tdive
ni}<v'<2u^{\tdive ni}$ for some $i<k$. If $\tdive ni=2\tdive n{i+1}$, we have $u'^2<v'<2u'^2$, hence P loses. Otherwise
$\tdive ni=2\tdive n{i+1}+1$; in a final round, C plays $u'^2$, and P responds with $u'^2\le u''<2u'^2$. Then
$u'^2<u''<2u'^2$ or $u'u''<v'<2u'u''$, thus P loses either way.

\ref{item:5}: In the first round, C makes P respond with a $v$ such that $u^n\le v<2u^n$. If $v=u^n$, C continues with
the strategy for $\powg(\vec u,u,u^n)$, otherwise with the strategy from~\ref{item:9}.
\end{Pf}
\begin{Rem}\label{rem:div-conq}
Extending Observation~\ref{obs:powg-div}, Lemma \ref{lem:div-conq}~\ref{item:5} implies the following for $n\ge2$: if $v\le u^n$
and $v\nmid u^n$, then $c(u,v)\le\cl{\log\fl{\log n}}+2$.
\end{Rem}
\begin{Def}\label{def:p-adic}
For any prime $p$ and $n\in\N_{>0}$, $\nu_p(n)$ denotes the $p$-adic valuation of~$n$: the maximal $k$ such that
$p^k\mid n$. (If it comes up, $\nu_p(0)$ is understood as $+\infty$.)
\end{Def}

Observe that any $n\in\N_{>0}\bez P_2^\N$ can be written uniquely as $n=2^{\nu_2(n)}v^r$, where $v$
(which is necessarily odd) is not a perfect power (which implies $v>1$), and $r>0$. We have
$r=\gcd\{\nu_p(n):\text{$p$ odd prime}\}$.
\begin{Thm}\label{thm:powg-ub}
Let $u=2^{\nu_2(u)}v^r$, where $v$ is not a perfect power, and let $d\nmid r$. Then
\begin{equation}\label{eq:27}
c(u)\le\cl{\log\fl{\log d}}+4.
\end{equation}
\end{Thm}
\begin{Pf}
In the first round, C can play $\fl{u^{1/d}}$ so that P responds with a $w$ such that $w\le u^{1/d}<2w$, hence
$2^iw^d<u<2^{i+1}w^d$ for some $i<d$: $2^iw^d=u$ is impossible as the odd part of $u$ is not a $d$th power. It remains to
show that for any such $w$, $c(u,w)\le\cl{\log\fl{\log d}}+3$.

Since $c(u,w^d,2^i)=0$, we have $c(u,w,2^i)\le\cl{\log\fl{\log d}}+2$ by Lemma \ref{lem:div-conq}~\ref{item:5}. One more
application of Lemma~\ref{lem:div-conq} gives
$c(u,w,2)\le\max\{\cl{\log\fl{\log d}}+3,\cl{\log\fl{\log i}}+2\}=\cl{\log\fl{\log d}}+3$, thus
$c(u,w)\le\cl{\log\fl{\log d}}+4$ (if C plays $2$, P has to respond with $2$).

We can improve this to $c(u,w)\le\cl{\log\fl{\log d}}+3$ by observing that $2$ is needed only in one branch. Mimicking
the proof of Lemma~\ref{lem:div-conq}, let C play $2^{i+1}-1$ so that P responds with a $z$ such that $2^i\le z<2^{i+1}$. If
$z=2^i$, C wins in $c(u,w,2^i)\le\cl{\log\fl{\log d}}+2$ more rounds for a total of $\cl{\log\fl{\log d}}+3$.
Otherwise, C makes P play $2$ in the second round, and wins in $\cl{\log\fl{\log i}}+1$ more rounds by
Lemma \ref{lem:div-conq}~\ref{item:9} for a total of $\cl{\log\fl{\log i}}+3$.
\end{Pf}
\begin{Rem}\label{rem:nu-nu}
Ignoring the exact constants, Theorem~\ref{thm:powg-ub} is equivalent to
\[c(u)\le\min\{\log\nu_q(\nu_p(u))+\log\log q:\text{$p$, $q$ primes, $p$ odd}\}+O(1).\]
\end{Rem}

Recall our convention that $\log x=0$ for $x<1$.

\begin{Thm}\label{cor:powg-ub}
Let $u=2^{\nu_2(u)}v^r$, where $v$ is not a perfect power, and $r>0$. Then
\begin{align}
\label{eq:28}c(u)&\le\cl{\log\log\log\log u}+4,\\
\label{eq:29}c(u)&\le\cl{\log\log\log r}+4,\\
\label{eq:30}c(u)&\le\cl{\log\log\log\nu_2(u)}+5.
\end{align}
\end{Thm}
\begin{Pf}
We start with~\eqref{eq:29}. There exists a $d\le n$ such that $d\nmid r$ if $r<\lcm\{1,\dots,n\}=e^{\psi(n)}$, where
$\psi(n)=\sum_{p^k\le n}\ln p$ is Chebyshev's function. By the prime number theorem, $\psi(n)\sim n$, hence we can find
$d\nmid r$ such that $d\le(1+o(1))\ln r$, thus $d\le\log r$ if $r$ is large enough. Then \eqref{eq:27} implies
\eqref{eq:29}.

To show that \eqref{eq:29} holds for all rather than just sufficiently large~$r$, we need to check small cases. First,
if $d=2,3$, then $\cl{\log\fl{\log d}}=0\le\cl{\log\log\log r}$, thus \eqref{eq:29} holds unless $6\mid r$, whence
$r\ge6>4$. Next, if $d\le7$, then $\cl{\log\fl{\log d}}=1\le\cl{\log\log\log r}$ (using $r>4$), thus \eqref{eq:29}
holds unless $2^2\times3\times5\times7=420\mid r$, whence $r\ge420>2^8$. Finally, it follows from known explicit bounds
on $\psi$ that $\lcm\{1,\dots,n\}>2^n$ for $n>8$ (in fact, it holds for $n\ge7$): Nagura~\cite{nagu:psi} proved
$\psi(n)>0{.}916n-2{.}318$ for all $n>0$, which implies $\psi(n)>n\ln2$ for $n\ge11$, and one can check the cases
$n=9,10$ by hand. Thus, if $r>2^8$, there is a $d\nmid r$ such that $d\le\cl{\log r}$, whence
$\cl{\log\fl{\log d}}\le\cl{\log\log\cl{\log r}}=\cl{\log\log\log r}$.

Since $r\le\log_3u\le\log u$, \eqref{eq:29} implies~\eqref{eq:28}.

For~\eqref{eq:30}, let C play $2^{\nu_2(u)+1}<u$ in the first round so that P responds with a $u'$ such that
$2^{\nu_2(u)}<u'\le2^{\nu_2(u)+1}$. If $u'=2^{\nu_2(u)+1}\nmid u$, then C can win in the next round by
Observation~\ref{obs:powg-div}. Otherwise, C can win in $\cl{\log\log\log\log_3u'}+4\le\cl{\log\log\log\nu_2(u)}+4$ further
rounds by (the proof of) \eqref{eq:28}.
\end{Pf}

\begin{Exm}\label{exm:c1}
$c(u)=1$ iff $u=5,6,7,17$.
\end{Exm}
\begin{Pf}
If $u=5,6,7$, C wins by playing $2$, forcing P to respond with $2$, as $2^2<u<2\cdot2^2$. If $u=17$, C plays $4$, and P
responds with $v=3,4$; then $v^2<u<2v^2$. Conversely, if $c(u)=1$, let $2x$ be the winning move of C (assumed even by
Remark~\ref{rem:powg-even}); then $v^2<u<2v^2$ or $u^2<v<2u^2$ for all $v\in(x,2x]$, i.e.,
$[x+1,2x]\sset\bigl[\fl{\sqrt{u/2}}+1,\cl{\sqrt u}-1\bigr]\cup[u^2+1,2u^2-1]$. There is a gap between the last two
intervals as $\cl{\sqrt u}-1<u<u^2+1$, thus $[x+1,2x]\sset\bigl[\fl{\sqrt{u/2}}+1,\cl{\sqrt u}-1\bigr]$ or
$[x+1,2x]\sset[u^2+1,2u^2-1]$. The latter makes $u^2\le x$ and $2x<2u^2$, which is impossible. The former amounts to
$\sqrt{u/2}-1<x<\frac12\sqrt u$; in particular, $\sqrt{2u}-2<\sqrt u$, thus $\sqrt u<2/(\sqrt2-1)=2(\sqrt2+1)$ and
$x<\sqrt2+1$, i.e., $x=1$ (in which case $4<u<8$) or $x=2$ (in which case $16<u<18$).
\end{Pf}

\begin{Exm}\label{exm:c2}
We have $c(u)\le2$ whenever $u$ satisfies one of the following conditions:
\begin{enumerate}
\item\label{item:6} $u>8$ and $16\nmid u$.
\item\label{item:10} The odd part of $u$ is not a square.
\item\label{item:11} $u<2304$ is not a power of $2$. (With some effort, one can check that $c(2304)=3$.)
\item\label{item:12} $u=\prod_{i<k}p_i^{e_i}$ for primes $p_0<\dots<p_{k-1}$, and there is $i<k$ such that
$p_i>2\prod_{j<i}p_j^{e_j}$.
\end{enumerate}
\end{Exm}
\begin{Pf}
Observe that C can force P to play $1$, $2$ (by playing the same), and in the first round, also $4$ (by playing $6$: if
P responds with $5,6$, C wins in the second round by Example~\ref{exm:c1}) and $8$ (by playing $8$; if P responds with
$5,6,7$, we use Example~\ref{exm:c1} again).

\ref{item:6}: C makes P play $8$ in the first round, and then plays $\cl{u/8}-1$, thus P responds with $v$ such that
$u/16<v<u/8$ ($v=u/16$ is impossible by assumption); C wins as $8v<u<2\cdot8v$.

\ref{item:10}: C plays $\fl{\sqrt u}$, thus P responds with $v$ such that $\frac12\sqrt u<v\le\sqrt u$; since $u\ne
v^2,2v^2$, we have $v^2<u<2v^2$ (and C wins) or $2v^2<u<4v^2$. In the latter case, C plays $\cl{u/v}-1$, thus P
responds with $w$ such that $w<u/v\le2w$. Then C wins as either $vw<u<2vw$, or $w=u/(2v)$ and $w^2<u<2w^2$.

\ref{item:11}: For $u=3$, C forces P to play $1$ and $2$. For $4<u<8$, C makes P play $2$ as in Example~\ref{exm:c1}. For
$8<u<64$, we can use \ref{item:6}, unless $u=48$, in which case we use~\ref{item:10}. For $64<u<128$, C makes P play
$8$ and wins as $8^2<u<2\cdot8^2$. For $128<u<256$, C makes P play $4$, and then plays $32$ so that P responds with
$16<v\le32$: either $4^2<v<2\cdot4^2$, or $v=32$ and $4\cdot32<u<2\cdot4\cdot32$, thus C wins. For $256<u<512$, C makes
P play $4$, and then plays $31$, thus P plays $16\le v<32$. Either $4^2<v<2\cdot4^2$, or $v=16$ and
$16^2<u<2\cdot16^2$. For $512<u<1024$, C makes P play $8$, and then plays $127$ so that either $8^2<v<2\cdot8^2$ or
$8\cdot64<u<2\cdot8\cdot64$. For $1024<u<2048$, C makes P play $4$, and then plays $32$ so that either
$4^2<v<2\cdot4^2$ or $32^2<u<2\cdot32^2$. (One can also do $4096<u<8192$ and $16384<u<32768$ using similar arguments.)
For $2048<u<2304=16\cdot144$, one of \ref{item:6} or~\ref{item:10} is applicable.

For $u=2304=48^2$, P can survive two rounds by playing in the first one an element of $\{u^n2^l:n\in\N,\abs
l\le4\}\cup\{u^{n+1/2}2^l:n\in\N,\abs l\le1\}$, but it is a bit tedious to check all cases.

\ref{item:12}: The assumption implies (and, actually, is equivalent to) that for some~$x<u$, namely
$x=\prod_{j<i}p_j^{e_j}$, there is no divisor $v\mid u$ such that $x<v\le2x$. Thus, C can play $2x$, and win in
the second round by Observation \ref{obs:powg-div}~\ref{item:7}.
\end{Pf}

The significance of point~\ref{item:12} of Example~\ref{exm:c2} is that the upper bounds from
Theorems \ref{thm:powg-ub} and~\ref{cor:powg-ub} cannot be asymptotically optimal: there are $u$ for which these bounds are arbitrarily
large, yet $c(u)=2$ (e.g., take $u=(2p)^{n!}$ for a large $n$, where $p>2^{n!+1}$ is a prime). Nevertheless, we will
show that the bounds are tight up to an additive constant under suitable conditions precluding \ref{item:12} and
similar cases (viz., in the decomposition $u=2^{\nu_2(u)}v^r$, $v$ is sufficiently smaller than $2^{\nu_2(u)}$).

We now come to the main technical part of our lower bound on $c(u)$. Recall that the $1$-norm of a vector
$\vec x\in\R^t$ is $\dlh{\vec x}_1=\sum_{i<t}\abs{x_i}$.
\begin{Lem}\label{lem:powg-lb}
Let $v\in\N_{>0}$, and define the sequences $\p{D_k,N_k,B_k:k\in\N_{>0}}$ by $D_1=1$, $N_1=3$, $B_1=0$,
$D_{k+1}=D_k\lcm\{1,\dots,N_k\}$, $N_{k+1}=N_k^2$, and $B_{k+1}=2N_kB_k+N_k^2\cl{D_k\log v}$. Then the following holds
for all $k\ge1$ and all $\vec u\in\N_{>0}^t$ of the form $u_i=2^{l_i}v^{r_i}$, $l_i,r_i\in\N$, for each $i<t$:

If
\begin{equation}\label{eq:31}
D_k\mid r_i
\end{equation}
for each $i<t$, and
\begin{equation}\label{eq:32}
\dlh{\vec n}_1\le N_k\et\sum_{i<t}n_ir_i>0\implies\sum_{i<t}n_il_i\ge B_k
\end{equation}
for all $\vec n\in\Z^t$, then $c(\vec u)\ge k$.
\end{Lem}
\begin{Pf}
We prove the statement by induction on $k$. We may assume $v$ is not a power of $2$ (whence $v\ge3$), as otherwise
$c(\vec u)=\infty$ trivially satisfies the conclusion.

For $k=1$, we have to show that there are no $h,i,j<t$ such that $u_iu_j<u_h<2u_iu_j$. Fixing $h,i,j$, put $\vec
n=e^h-e^i-e^j$, where $e^g$ denotes the $g$th standard unit vector (i.e., if $h,i,j$ are distinct, then $n_h=1$ and
$n_i=n_j=-1$). Clearly, $\dlh{\vec n}_1\le3=N_1$, thus we may apply~\eqref{eq:32}: if $r_h-r_i-r_j>0$, then
$l_h-l_i-l_j\ge0=B_1$, hence $u_h/(u_iu_j)=2^{l_h-l_i-l_j}v^{r_h-r_i-r_j}\ge v>2$. Likewise, if $r_h-r_i-r_j<0$, we
may apply \eqref{eq:32} to $-\vec n$, and obtain $u_h/(u_iu_j)\le v^{-1}<1$. Finally, if $r_h-r_i-r_j=0$, then
$u_h/(u_iu_j)$ is a power of~$2$, hence it cannot be strictly between $1$ and~$2$.

Assume the statement holds for $k$, and that $\vec u$, $\vec l$, and $\vec r$ satisfy \eqref{eq:31} and~\eqref{eq:32}
for $k+1$ in place of~$k$. Using the induction hypothesis, it suffices to show that for every $x\ge1$, there exists
$u_t=2^{l_t}v^{r_t}$ such that $u_t\le x<2u_t$, and $\p{\vec u,u_t}\in\N_{>0}^{t+1}$ satisfies \eqref{eq:31} and
\eqref{eq:32} for~$k$. We will write $r=r_t$ and $l=l_t$ for short. Observe that $u_t\le x<2u_t$ amounts to
$l=\fl{\log x-r\log v}$, thus we can only vary~$r$; we will check that \eqref{eq:31} and~\eqref{eq:32} translate to
conditions on $r\in\Z$ that are satisfiable together, using our assumptions on~$\vec u$. (We also need to ensure
$r,l\ge0$, but this easily follows from~\eqref{eq:32}, hence we need not worry about it.)

Condition~\eqref{eq:31} clearly holds for $i<t$ as $D_k\mid D_{k+1}$, thus we only need to make sure $r$ is a multiple
of $D_k$. Condition~\eqref{eq:32} also holds automatically when $n_t=0$, as $N_k\le N_{k+1}$ and $B_k\le B_{k+1}$.
The other cases give upper or lower bounds on~$r$, depending on the sign of $n_t$. For $n_t>0$ (renamed to~$n$, and the
rest of $\vec n$ negated), \eqref{eq:32} amounts to
\[\dlh{\vec n}_1+n\le N_k\et r>\sum_{i<t}\frac{n_ir_i}n\implies l\ge\frac{B_k}n+\sum_{i<t}\frac{n_il_i}n\]
for all $\vec n\in\Z^t$ and $n\in\N_{>0}$, i.e.,
\[\dlh{\vec n}_1+n\le N_k\implies r\le\sum_{i<t}\frac{n_ir_i}n\quad\text{or}\quad\log x-r\log
v\ge\CL{\frac{B_k}n+\sum_{i<t}\frac{n_il_i}n}.\]
The largest integer multiple $r=D_kr'$ that satisfies this condition is characterized by
\[r'\le U_{\vec n,n}:=\max\left\{\sum_{i<t}\frac{n_ir_i}{D_kn},
  \Fl{\frac1{D_k\log v}\left(\log x-\CL{\frac{B_k}n+\sum_{i<t}\frac{n_il_i}n}\right)}\right\}\]
for all $\vec n\in\Z^t$ and $n\in\N_{>0}$ such that $\dlh{\vec n}_1+n\le N_k$, using the fact that for $1\le n\le N_k$,
$D_kn\mid D_{k+1}\mid r_i$.

Likewise, the cases of~\eqref{eq:32} with $n_t<0$ (renamed to $-m$, and the rest of $\vec n$ to $\vec m$) amount to
\[\dlh{\vec m}_1+m\le N_k\implies r\ge\sum_{i<t}\frac{m_ir_i}m\quad\text{or}\quad\log x-r\log
v<\Fl{-\frac{B_k}m+\sum_{i<t}\frac{m_il_i}m}+1\]
for all $\vec m\in\Z^t$ and $m\in\N_{>0}$, and the least multiple of $D_k$ with this property is characterized by
\[r'\ge L_{\vec m,m}:=\min\left\{\sum_{i<t}\frac{m_ir_i}{D_km},
  \Fl{\frac1{D_k\log v}\left(\log x+\CL{\frac{B_k}m-\sum_{i<t}\frac{m_il_i}m}-1\right)}+1\right\}\]
for all $\vec m\in\Z^t$ and $m\in\N_{>0}$ such that $\dlh{\vec m}_1+m\le N_k$. Thus, an $r$ that satisfies all the
necessary conditions exists iff for every $\vec n,\vec m\in\Z^t$ and $n,m\in\N_{>0}$,
\begin{equation}\label{eq:36}
\dlh{\vec n}_1+n\le N_k\et\dlh{\vec m}_1+m\le N_k\implies L_{\vec m,m}\le U_{\vec n,n}.
\end{equation}
This clearly holds if $\frac1m\sum_im_ir_i\le\frac1n\sum_in_ir_i$, hence we may assume
$\frac1m\sum_im_ir_i>\frac1n\sum_in_ir_i$. Then the assumption \eqref{eq:32} for $\vec u$, applied to $n\vec m-m\vec
n$, implies
\[\sum_{i<t}\frac{m_il_i}m-\sum_{i<t}\frac{n_il_i}n\ge\frac{B_{k+1}}{nm}\ge\frac{B_k}m+\frac{B_k}n+\cl{D_k\log v},\]
using the bounds
\[\dlh{n\vec m-m\vec n}_1\le n\dlh{\vec m}_1+m\dlh{\vec n}_1
  \le\bigl(n+\dlh{\vec n}_1\bigr)\bigl(m+\dlh{\vec m}_1\bigr)\le N_k^2=N_{k+1}\]
and
\[B_{k+1}=2N_kB_k+N_k^2\cl{D_k\log v}\ge(n+m)B_k+nm\cl{D_k\log v}.\]
It follows that
\[\CL{\frac{B_k}n+\sum_{i<t}\frac{n_il_i}n+\frac{B_k}m-\sum_{i<t}\frac{m_il_i}m}\le-\cl{D_k\log v}\le-D_k\log v,\]
thus
\[\CL{\frac{B_k}n+\sum_{i<t}\frac{n_il_i}n}+\CL{\frac{B_k}m-\sum_{i<t}\frac{m_il_i}m}\le1-D_k\log v\]
and
\[\left(\log x+\CL{\frac{B_k}m-\sum_{i<t}\frac{m_il_i}m}-1\right)+D_k\log v
  \le\log x-\CL{\frac{B_k}n+\sum_{i<t}\frac{n_il_i}n}.\]
This yields~\eqref{eq:36}.
\end{Pf}

\begin{Thm}\label{thm:powg-lb}
Let $u=2^lv^r$, where $v>1$, $r>0$, and $l/\log v\ge10^8$. Then
\begin{equation}\label{eq:33}
c(u)\ge\min\bigl\{\fl{\log\cl{\log_3d}}+1:d\nmid r\bigr\}\cup\bigl\{\fl{\log\log_3\log_4(l/\log v)}+3\bigr\}.
\end{equation}
We may use the simpler bounds $\cl{\log\log_3d}$ or $\fl{\log\cl{\log d}}$ in place of $\fl{\log\cl{\log_3d}}+1$.
\end{Thm}
\begin{Pf}
Applying Lemma~\ref{lem:powg-lb} with $t=1$, we see that $c(u)\ge k$ whenever $D_k\mid r$ and $l\ge B_k$; it remains to
estimate these quantities. Expanding the recurrences, we have
\[N_k=3^{2^{k-1}},\qquad D_k=\prod_{i=0}^{k-2}L\bigl(3^{2^i}\bigr),\qquad\text{where }L(n)=\lcm\{1,\dots,n\}.\]
Observe that $L(n)L(m)\mid L(nm)$: whenever $1\le a\le n$ and $1\le b\le m$, we have $ab\mid L(nm)$. Thus,
\[D_k\Bigm|L\Bigl(\prod_{i=0}^{k-2}3^{2^i}\Bigr)=L\bigl(3^{2^{k-1}-1}\bigr),\]
and a sufficient condition for $D_k\mid r$ is that $d>3^{2^{k-1}-1}$ for all $d\nmid r$. Since
\[d>3^{2^{k-1}-1}\iff\cl{\log_3d}\ge2^{k-1}\iff\log\cl{\log_3d}+1\ge k,\]
we see that $D_k\mid r$ holds for any $k$ such that
\[k\le\min\{\fl{\log\cl{\log_3d}}+1:d\nmid r\}.\]
We also observe that $\fl{\log\cl{\log_3d}}+1\ge\cl{\log\cl{\log_3d}}=\cl{\log\log_3d}$ and
$\fl{\log\cl{\log_3d}}+1=\fl{\log(2\cl{\log_3d})}\ge\fl{\log\cl{\log d}}$, as $\log d\le2\log_3d$ implies $\cl{\log
d}\le2\cl{\log_3 d}$. 

The recurrence for $B_k$ resolves to
\[B_k=N_k\sum_{i=1}^{k-1}2^{k-1-i}\cl{D_i\log v}\le3^{2^{k-1}}\sum_{i=1}^{k-1}2^{k-1-i}(D_i+1)\log v=:B'_k\log v.\]
Recalling Chebyshev's function from the proof of Theorem~\ref{cor:powg-ub}, we have
$\ln L\bigl(3^{2^j}\bigr)=\psi\bigl(3^{2^j}\bigr)\sim3^{2^j}$, whence
$\ln D_i=\sum_{j\le i-2}\psi\bigl(3^{2^j}\bigr)\sim3^{2^{i-2}}$. It follows that the above sum for~$B'_k$ is dominated
by the $i=k-1$ term, and $\ln B'_k\sim3^{2^{k-3}}$; thus, for all sufficiently large $k$, $\log_4B'_k<3^{2^{k-3}}$ and
$k>3+\log\log_3\log_4B'_k$. That is, $l\ge B_k$ if $k\le\fl{\log\log_3\log_4(l/\log v)}+3$, provided the latter is
large enough.

For an explicit bound, we claim that $\log\log_3\log_4B'_k<k-3$ for all $k\ge5$, thus \eqref{eq:33} holds whenever
$l/\log v\ge B'_4=99{,}353{,}223$. For $k=5$, direct computation gives
$B'_5\approx6{.}333\times10^{46}$ and $\log\log_3\log_4B'_5<1{.}9865$ (this can be verified e.g.\ in Sage, along with the
value of $B'_4$ given above). Assume $k\ge6$. Nagura~\cite{nagu:psi} proved $\psi(n)<1{.}086n$ for all $n>0$, thus
\[\log_4D_{k-1}<1{.}086\log_4e\sum_{i\le k-3}3^{2^i}<0{.}7834\,\alpha_{k-3}\,3^{2^{k-3}},\]
where $\alpha_j=\sum_{i\le j}3^{2^i-2^j}$. For $j\ge1$, $3^{2^j}(\alpha_j-1)=\sum_{i<j}3^{2^i}\ge3>\alpha_j$, thus
$\alpha_{j+1}=1+3^{-2^j}\alpha_j<\alpha_j$, i.e., $\alpha_j$ is decreasing. Since
$\alpha_3=1+3^{-4}+3^{-6}+3^{-7}=2218/2187<1{.}0142$, we get
\[\log_4D_{k-1}<0{.}7834\,\alpha_3\,3^{2^{k-3}}<0{.}795\times3^{2^{k-3}}\]
and
\[\log_4(N_kD_{k-1})<0{.}795\times3^{2^{k-3}}+2^{k-1}\log_43\le\left(0{.}795+\frac{2^5}{3^{2^3}}\log_43\right)3^{2^{k-3}}
  <0{.}799\times3^{2^{k-3}}.\]
Since $D_{i+1}/D_i=L\bigl(3^{2^{i-1}}\bigr)\ge6$, we have
\begin{align*}
\frac{B'_k}{N_kD_{k-1}}-1&=\frac{2^{k-1}-1}{D_{k-1}}+\sum_{i=1}^{k-2}2^{k-1-i}\frac{D_i}{D_{k-1}}<
  \frac{2^{k-1}-1}{D_{k-1}}+3\frac{D_{k-2}}{D_{k-1}}\\
&\le\frac{31}{D_5}+\frac3{L\bigl(3^{2^3}\bigr)}<4\times10^{-2846},
\end{align*}
thus $\log_4(B'_k/(N_kD_{k-1}))<4\times10^{-2846}\log_4e<3\times10^{-2846}$, and $\log_4B'_k<0{.}8\times3^{2^{k-3}}$.
\end{Pf}

\begin{Exm}\label{exm:6-2-2-k-fac}
$k+1\le c\bigl(6^{2^{\scriptstyle2^k}\textstyle!}\bigr)\le k+4$ for all $k\ge0$.
\end{Exm}
\begin{Pf}
We write $u=6^{2^{\scriptstyle2^k}\textstyle!}=2^lv^r$ where $v=3$ and $l=r=2^{2^k}!$. The least $d$ not dividing
$r$ is the least prime larger than $2^{2^k}$, thus $2^{2^k}<d<2^{2^k+1}$ by Bertrand's postulate, and $c(u)\le k+4$ by
Theorem~\ref{thm:powg-ub}. For a lower bound, we have $c(36)=2$ and $c(6^4)\ge2$ by Examples \ref{exm:c1} and~\ref{exm:c2}, thus we may assume
$k\ge2$. Since $n!>4^n$ for $n\ge9$, we have
$4^{3^{2^{k-1}}}\log3<4^{3^{2^{k-1}}+1}<4^{4^{2^{k-1}}}=4^{2^{2^k}}<2^{2^k}!$ for $k\ge2$, i.e.,
$\log\log_3\log_4(2^{2^k}!/\log3)>k-1$, and $2^{2^k}/\log3>10^8$. Thus, Theorem~\ref{thm:powg-lb} gives
$c(u)\ge\fl{\log\cl{\log_3d}}+1\ge k$.

To improve this to $c(u)\ge k+1$, we may apply Lemma~\ref{lem:powg-lb} directly. We know $l\ge B_{k+2}$ from the proof of
Theorem~\ref{thm:powg-lb}, thus it suffices to show that $D_{k+1}\mid r$. Any prime $p\mid D_{k+1}$ is bounded by
$3^{2^{k-1}}$; since $\nu_p(L(n))=\fl{\log_pn}$, we have
\[\nu_p(D_{k+1})=\sum_{i=0}^{k-1}\Fl{\log_p3^{2^i}}\le\Fl{\sum_{i=0}^{k-1}2^i\frac{\log3}{\log p}}\le\Fl{\frac{2^k}{\log_3p}},\]
while $\nu_p(r)=\sum_{i\ge1}\fl{2^{2^k}/p^i}\ge\fl{2^{2^k}/p}$. It remains to observe that $2^{2^k}/p\ge2^k/\log_3p$
as $p/\log_3p\le3^{2^{k-1}}/2^{k-1}\le2^{2^k}/2^k$.
\end{Pf}

\begin{Thm}\label{thm:hier}
If $T$ is any $\Sigma_1$-sound $\langor$-theory, then the formulas $\{\pwin^1_k(u):k\ge1\}$ are pairwise inequivalent
over~$T$.
\end{Thm}
\begin{Pf}
Since $T$ remains $\Sigma_1$-sound after adding any set of true $\Pi_1$ sentences, we may assume $T\Sset\idz$. In view
of Lemma~\ref{lem:theta}, it suffices to prove $T\nvdash\pwin^1_k(u)\to\pwin^1_{k+1}(u)$ for any $k\ge1$. By
Example~\ref{exm:6-2-2-k-fac} and Lemma~\ref{lem:c-init}, there exists $n\in\N$ such that $c(n)=k$, i.e.,
$\N\model\pwin^1_k(\ob n)\land\neg\pwin^1_{k+1}(\ob n)$. We observe that the existential quantifiers in
Definition~\ref{def:teip} can be bounded with $u_i\le x_i$, thus $\neg\pwin^1_{k+1}$ is equivalent to a $\Sigma$-formula (i.e.,
a formula built using existential and bounded universal quantifiers from a $\Delta_0$ formula). It follows that
$\neg\pwin^1_{k+1}(\ob n)$ is provable in $\thry Q\sset T$, and if we assume that $T$ is $\Sigma$-sound, then
$T\nvdash\neg\pwin^1_k(\ob n)$.

If we only have the weaker assumption that $T$ is $\Sigma_1$-sound, we need to be a bit more careful. If we fix
$t,k\ge1$ and $v>1$, then $D_k$, $N_k$, and $B_k$ are constants, and properties \eqref{eq:31} and~\eqref{eq:32} can be
written as a $\Delta_0$~formula $\alpha^t_{v,k}(\vec u)$: $\vec r$ and~$\vec l$ are bounded by~$\vec u$, and we can
express the condition $u_i=2^{l_i}v^{r_i}$ by a bounded formula as the graph of powering $x^y=z$ is
$\Delta_0$-definable \cite[\S V.3(c)]{hp}; also, there is only a constant number of choices for~$\vec n$. Then the
$\Pi_1$ sentences
\begin{align}
\label{eq:34}&\forall\vec u\:\bigl(\alpha^t_{v,1}(\vec u)\to\pwin^t_t(\vec u)\bigr),\\
\label{eq:35}&\forall\vec u,x\:
  \bigl(\alpha^t_{v,k+1}(\vec u)\land x>0\to\exists u_t\le x\:(x<2u_t\land\alpha^{t+1}_{v,k}(\vec u,u_t))\bigr)
\end{align}
are true in~$\N$ and imply $\alpha^t_{v,k}(\vec u)\to\pwin^t_{t+k-1}(\vec u)$. Thus, if we fix $k$ and
$n=6^{2^{\scriptstyle2^{k-1}}\textstyle!}$ from Example~\ref{exm:6-2-2-k-fac}, there is a true $\Pi_1$ sentence that
implies $\pwin^1_k(\ob n)$, while $\thry Q\vdash\neg\pwin^1_{k+4}(\ob n)$. Using the $\Sigma_1$-soundness
of~$T$, there is a model $\sM\model T+\pwin^1_k(\ob n)+\neg\pwin^1_{k+4}(\ob n)$. By the argument in
Lemma~\ref{lem:c-init}, there is $m\le n$ such that $\sM\model\pwin^1_k(\ob m)\land\neg\pwin^1_{k+1}(\ob m)$.
\end{Pf}
\begin{Rem}\label{rem:vtc0}
Without going into the details, we claim that the lower bound in Lemma~\ref{lem:powg-lb} can be formalized for standard
$k$, $t$, and $v$ in the theories $\idz$ and $\dicr$ (or equivalently, $\vtc$); that is, these theories prove
\eqref{eq:34} and~\eqref{eq:35}. It follows that Theorem~\ref{thm:hier} holds for all consistent extensions of $\idz$ or
$\dicr$, regardless of their $\Sigma_1$-soundness.
\end{Rem}

\section{Oddless interpretation}\label{sec:oddless}

We determined the $\langor$-fragment of $\teipp$ to be $\teip$ in Section~\ref{sec:teip}. But for completeness, we
mention that there is another natural approach of relating $\teipp$ to $\langor$-theories which places an upper bound
on the strength of $\teip$: it is common to define the set of powers of~$2$ in arithmetical theories by the formula
\[\powii(u)\iff\forall x\:(x\mid u\to x=1\lor2\mid x),\]
expressing that $u$ has no nontrivial odd divisors (hence we may call such elements ``oddless''). Let $\pi_2$ denote
the translation of $\langp$ into $\langor$ which is absolute on $\langor$, and interprets $P_2$ by $\powii$: i.e., for
any $(\langp)$-formula $\fii(\vec x)$, let $\fii^{\pi_2}(\vec x)$ denote the $\langor$-formula obtained by replacing
each subformula of the form $P_2(t)$ with $\powii(t)$; an $\langor$-theory $T$ interprets a $(\langp)$-theory $S$
via~$\pi_2$ if $T\vdash\fii^{\pi_2}$ for each axiom (and therefore, each theorem) $\fii$ of~$S$. We are particularly
interested in relating $\teip$ and friends to standard fragments of bounded arithmetic using~$\pi_2$.

Recall that $x\mid y$ is an $E_1$~formula equivalent to an $U_1$~formula over~$\io$, thus $\powii$ is equivalent to a
$U_1$~formula (the quantifier over~$x$ in the definition can be clearly bounded by~$u$).
\begin{Obs}\label{obs:pow2}
$\io$ proves $\powii(u)\land v\mid u\to\powii(v)$.
\noproof\end{Obs}
\begin{Thm}\label{thm:oddless}
The smallest theory that interprets $\teipp$ via~$\pi_2$ is the theory $\teippow$, extending $\io$ by the axioms
\begin{align}
\label{eq:9}\tag{$\powii$-IP}&\forall x>0\:\exists u\:\bigl(\powii(u)\land u\le x<2u\bigr),\\
\label{eq:10}\tag{$\powii$-Div}&\forall u,v\:\bigl(\powii(u)\land\powii(v)\land u\le v\to u\mid v\bigr).
\end{align}
It is included in $\thry{IE_1}+{}$
\begin{equation}\label{eq:37}\tag{$\powii$-Cof}
\forall x\:\exists u>x\:\powii(u),
\end{equation}
in~$\thry{IE_2}$, and in $\dicr$ \textup(or equivalently, $\vtc$\textup).
\end{Thm}
\begin{Pf}
\eqref{eq:9} and~\eqref{eq:10} are almost literally the $\pi_2$-translations of axioms \eqref{eq:4}
and~\eqref{eq:5}. The only difference is that to
get~\eqref{eq:5}, we should also require $\powii(v/u)$; but this follows from Observation~\ref{obs:pow2}.

It is well known that $\thry{IE_1}$ proves that any two integers have a gcd (even with B\'ezout cofactors; see the
argument in Wilmers~\cite[Lemma~2.4]{wilmers}). This implies \eqref{eq:10}: if $\powii(u)$ and $\powii(v)$, let
$d=\gcd(u,v)$. Then $u'=u/d$ and $v'=v/d$ are coprime, hence one of them is odd. Being divisors of the oddless $u$
or~$v$, this implies $u'=1$ or $v'=1$, i.e., $u\mid v$ or $v\mid u$.

Working in $\thry{IE_1}$, let $x>0$ be given, and assume there exists an oddless $v>x$.
The $E_1$~formula $\fii(u)\equiv\exists u'\le x\,(u'\ge u\land u'\mid v)$ satisfies $\fii(0)\land\neg\fii(x+1)$, thus
using $E_1$-induction, there is $u$ such that $\fii(u)\land\neg\fii(u+1)$; then $u$ is the largest divisor of~$v$ such
that $u\le x$. Since $v/u>1$ must be even, we have $2u\mid v$, hence $x<2u$ by the maximality of~$u$, and $\powii(u)$ by
Observation~\ref{obs:pow2}. Thus, $\thry{IE_1}+\eqref{eq:37}\vdash\eqref{eq:9}$.

$\thry{IE_2}$ proves the $E_2$~formula $\exists u\le2x\,(u>x\land\powii(u))$ by induction on~$x$, as it is easy to see
that $\powii(u)$ implies $\powii(2u)$.

In $\vtc$, there is a canonical $2^n$ function from unary to binary integers, and every binary $X>0$ can be written as
$X=2^nX'$ with $X'$ odd. It follows easily that $\powii(X)$ holds iff $X$ is in the image of~$2^n$. Then
\eqref{eq:9} follows as $2^{n-1}\le X<2^n$ for $n$ given by the length function $|X|$ of $\vtc$, and
\eqref{eq:10} follows from $2^m=2^n2^{m-n}$ (for $n\le m$).
\end{Pf}

\begin{Cor}\label{cor:oddless-teip}
The theories $\teippow$, $\thry{IE_1}+\eqref{eq:37}$, $\thry{IE_2}$, and $\dicr$ contain $\teip$.

Any model of any of these theories has an elementary extension $\sM=\p{M,\dots}$ which is an
EIP of a RCEF $\p{\sR,\exp}$ satisfying GA such that $\exp[M]=\{u\in M:\sM\model\powii(u)\}$.
\noproof\end{Cor}

We mention that Corollary~\ref{cor:oddless-teip} does not quite reprove the main result of \cite{ej:vtceip}, which guarantees
that every countable model of $\dicr$ is outright an EIP of a RCEF satisfying GA, without taking an elementary
extension first.

If $T$ is $\teippow$ or any of the stronger $\langor$-theories from Theorem~\ref{thm:oddless}, and $\sM\model T$, it is a
natural question whether the expansion $\p{\sM,\powii^\sM}$ of $\sM$ to a model of $\teipp$ is unique: is it necessary
that $P_2$ consists of oddless numbers for \emph{every} expansion of a model of $T$ to a model of $\teipp$, perhaps if
$T$ is sufficiently strong? In other words, is $T+\teipp$ equivalent to the expansion of $T$ by the definition
$P_2(u)\eq\powii(u)$?

Our results from the previous section give a negative answer, even when $T$ is as strong as the true arithmetic:
\begin{Thm}\label{thm:nonunique}
Let $T$ be a $\Sigma_1$-sound $\langor$-theory. Then
\[T+\teipp\nvdash P_2(u)\to3\nmid u,\]
thus there exists a model $\p{\sM,P_2^\sM}\model T+\teipp$ such that $P_2^\sM\nsset\powii^\sM\nsset P_2^\sM$.
\end{Thm}
\begin{Pf}
Since adding true $\Pi_1$ sentences preserves $\Sigma_1$-soundness, we may assume $T\Sset\idz\Sset\teippow$.
By the proof of Theorem~\ref{thm:hier}, for each $k$ there exists $n\in\N$ divisible by $3$ (in fact, a power of $6$) such that
$T+\pwin^1_k(\ob n)$ is consistent. Thus by compactness, there exists a countable $\sM\model T$ and $u\in M$ such that
$3\mid u$ and $\sM\model\{\pwin^1_k(u):k\ge1\}$. We may assume $\sM$ to be recursively saturated; then by
the proof of Theorem~\ref{thm:teip}, there exists $P_2\sset M$ such that $\p{\sM,P_2}\model\teipp$ and $u\in P_2$. Clearly,
$\sM\nmodel\powii(u)$. On the other hand, since $\p{\sM,\powii^\sM}\model\teipp$, there is $v<u<2v$ such that
$\sM\model\powii(v)$, and we cannot have $v\in P_2$.
\end{Pf}

Let us mention that we do get uniqueness (for sufficiently strong $T$, namely including $\eqref{eq:10}$) if we
extend $\teipp$ with an axiom ensuring the downward closure of $P_2$ under divisibility:
\begin{Prop}\label{prop:unique}
An $\langor$-structure $\sM\model\eqref{eq:10}$ has at most one expansion to a model of $\teipp+{}$
\begin{equation}\label{eq:38}\tag{$P_2$-Down}
\forall u,v\:\bigl(P_2(u)\land v\mid u\to P_2(v)\bigr),
\end{equation}
namely $\p{\sM,\powii^\sM}$.
\end{Prop}

This follows from the next lemma, characterizing extensions of $\teipp$ in which $P_2$ is provably defined
by~$\powii$.
\begin{Lem}\label{lem:pow2-def}
These theories are equivalent:
\begin{enumerate}
\item\label{item:13} $\teipp+\forall u\,\bigl(P_2(u)\eq\powii(u)\bigr)$.
\item\label{item:14} $\teipp+\forall u\,\bigl(P_2(u)\to\powii(u)\bigr)+\eqref{eq:10}$.
\item\label{item:15} $\teipp+\forall u\,\bigl(\powii(u)\to P_2(u)\bigr)+\eqref{eq:37}$.
\item\label{item:16} $\teipp+\eqref{eq:38}+\eqref{eq:10}$.
\end{enumerate}
\end{Lem}
\begin{Pf}

\ref{item:13}\txto\ref{item:16}: \eqref{eq:38} follows from $\forall u\,\bigl(P_2(u)\eq\powii(u)\bigr)$ using
Observation~\ref{obs:pow2}, and $\forall u\,\bigl(\powii(u)\to P_2(u)\bigr)$ and \eqref{eq:5} imply
\eqref{eq:10}.

\ref{item:16}\txto\ref{item:14}: We will show $\teipp+\eqref{eq:38}\vdash P_2(u)\to\powii(u)$. Assume $P_2(u)$,
and let $v\mid u$ be such that $v>1$. Then $P_2(v)$ by \eqref{eq:38}, and $P_2(2)$ by Lemma~\ref{lem:teipp}, thus $v$
is even by \eqref{eq:5}.

\ref{item:14}\txto\ref{item:13}: We need to show $\powii(u)\to P_2(u)$. Assuming $\powii(u)$,
\eqref{eq:4} gives $v\le u<2v$ such that $P_2(v)$. Then $\forall u\,\bigl(P_2(u)\to\powii(u)\bigr)$ yields $\powii(v)$,
thus $v\mid u$ by \eqref{eq:10}. Since $1\le u/v<2$, we obtain $v=u$, i.e., $P_2(u)$.

\ref{item:13}\txto\ref{item:15}: $\teipp$ proves $\forall x\,\exists u>x\,P_2(u)$, which together with $\forall
u\,\bigl(P_2(u)\to\powii(u)\bigr)$ yields \eqref{eq:37}.

\ref{item:15}\txto\ref{item:13}: We need to show $P_2(u)\to\powii(u)$. Assuming $P_2(u)$, let $v>u$ be such that
$\powii(v)$. Then $P_2(v)$ as well, hence $u\mid v$ by \eqref{eq:5}. We get $\powii(u)$ by Observation~\ref{obs:pow2}.
\end{Pf}

\section{Discussion}\label{sec:conclusion}

We managed to axiomatize the first-order consequences of being an EIP of RCEF\@. While we obtained a simple finite list
of ``obvious'' axioms in languages including $2^x$ of~$P_2$, the axiomatization in the basic language of arithmetic
involves an unexpected infinite schema of axioms expressing the existence of winning strategies in~$\powg$.

Unlike the original Shepherdson's theorem, our results only characterize EIP of RCEF up to elementary equivalence, as
models of the resulting first-order theories may require an elementary extension to become an EIP of a RCEF. (We know
this is sometimes necessary for models of $\teip$ and $\teipp$ from Example~\ref{exm:shep}; we do not have a similar example
for $\teipe$, but it seems very likely that it should exist as well.) We leave open the problem whether a more precise
characterization is possible, at least for countable structures.

Another problem we left open is whether $\teip$ is finitely axiomatizable over~$\io$. Theorem~\ref{thm:hier} provides some
heuristic support for a negative answer, though the evidence it provides is quite limited (notice that
Theorem~\ref{thm:hier} exhibits a strict hierarchy even over $\Th(\N)$, whereas $\teip$ clearly \emph{is} finitely
axiomatizable over sufficiently strong theories, e.g.\ $\thry{IE_2}$).

While $\teip$ is a strict extension of~$\io$, it is not quite clear how much stronger it really is. In terms of literal
inclusion of theories, $\teip$ is contained in $\thry{IE_2}$ and $\dicr$, but we do not know if it is contained in
$\thry{IE_1}$. But perhaps a better assessment of the relative strength of $\teip$ is to estimate the minimal
complexity of sentences separating $\teip$ from~$\io$. In particular, a problem suggested by L.~Ko\l odziejczyk is to
determine what Diophantine equations are solvable in (extensions with negatives of) models of $\teip$, and whether they
are the same as those solvable in models of $\io$ (or equivalently, in $\Z$-rings, cf.\ Wilkie~\cite{wil:iop}); recall
that it is an old open question, going back to Shepherdson~\cite{sheph}, whether solvability of Diophantine equations
in models of $\io$ is decidable. A closely related question is whether $\teip$ is $\forall_1$-conservative over $\io$.

\providecommand\gobble[1]{}
\ifx\url\undefined
  {\catcode`\/=13 \gdef/{\string/\futurelet\nexttoken\finishslash}
  \gdef\finishslash{\ifx\nexttoken/\else\penalty\relpenalty\fi}}
  \def\url{\begingroup\catcode`\~=12 \catcode`\_=12 \catcode`\/=13 \finishurl}
  \def\finishurl#1{\texttt{#1}\endgroup}
\fi
\providecommand\bysame{\leavevmode\hbox to5em{\hrulefill}\thinspace}

\bibliographystyle{mybib}
\bibliography{teip}

\end{document}